\documentclass[11pt]{article}
\usepackage[margin=1in]{geometry}
\usepackage{amsmath,amsfonts,amssymb,graphics,graphicx,amsthm}
\usepackage{color}
\usepackage{hyperref}
\hypersetup{
    colorlinks=true,
    linkcolor=blue,
    citecolor=red,
    urlcolor=blue,
    pdfborder={0 0 0}
}

\newtheorem{theorem}{Theorem}[section]
\newtheorem{corollary}[theorem]{Corollary}
\newtheorem{lemma}[theorem]{Lemma}
\newtheorem{proposition}[theorem]{Proposition}

\theoremstyle{remark}
\theoremstyle{remark}\newtheorem{remark}[theorem]{Remark}

\def\measure{\mu_{\sigma,h}}

\newcommand{\rI}{\mathrm{I}}
\newcommand{\rII}{\mathrm{II}}
\newcommand{\rIII}{\mathrm{III}}

\newcommand{\C}{\mathbb{C}}
\newcommand{\D}{\mathbb{D}}

\newcommand{\E}{\mathbb{E}}

\newcommand{\R}{\mathbb{R}}
\renewcommand{\P}{\mathbb{P}}

\newcommand{\eps}{\varepsilon}

\newcommand{\1}{\mathbf{1}}

\newcommand\MR[1]{\relax\ifhmode\unskip\spacefactor3000 \space\fi
	\MRhref{\expandafter\@rst #1 other}{#1}}
\newcommand{\MRhref}[2]{\href{http://www.ams.org/mathscinet-getitem?mr=#1}{MR#2}}

\DeclareMathOperator{\dist}{dist}

\DeclareMathOperator{\supp}{supp}
\DeclareMathOperator{\har}{har}
\DeclareMathOperator{\GFF}{GFF}
\DeclareMathOperator{\Cov}{Cov}
\DeclareMathOperator{\Var}{Var}

\def\fh{\mathfrak{h}}

\def\cE{\mathcal{E}}

\def\cB{\mathcal{B}}
\def\cA{\mathcal{A}}

\usepackage{fullpage}
\usepackage{comment}

\def\P{\mathbb{P}}
\def\E{\mathbb{E}}
\def\eps{\varepsilon}

\DeclareMathOperator{\var}{Var}
\DeclareMathOperator{\cov}{Cov}

\DeclareMathOperator{\SLE}{SLE}

\newcommand{\indic}[1]{\mathbf{1}_{\{#1\}}}
\parskip \medskipamount

\def\D{\mathbb{D}}

\begin {document}
\author{\begin{tabular}{c}Nathana\"el Berestycki\\[-4pt]\small University of Vienna \end{tabular} \and \begin{tabular}{c}Scott Sheffield\\[-4pt]\small MIT \end{tabular}\and \begin{tabular}{c}Xin Sun\\[-4pt]\small University of Pennsylvania\end{tabular}
}

\title{ Equivalence of Liouville measure and Gaussian free field}
\date{\today}
\maketitle

\begin{abstract}
	Given an instance $h$ of the Gaussian free field on a planar domain
	$D$ and a constant $\gamma \in (0,2)$, one can use various
	regularization procedures to make sense of the {\em Liouville quantum
		gravity area measure} $\mu := e^{\gamma h(z)} dz.$
	It is known that the field $h$ a.s.\ determines the measure $\mu_h$.
	We show that the converse is true: namely, $h$ is measurably determined by $\mu_h$.
	More generally, given a random closed fractal subset $\cA$ endowed with a Frostman measure {$\sigma$} whose support is $\cA$ (independent of $h$), a Gaussian multiplicative chaos measure {$\measure$} can be constructed.
	We give a mild condition on {$(\cA,\sigma)$} under which {$\measure$} determines $h$ restricted to $\cA$, in the sense that it determines its harmonic extension off $\cA$.
	Our condition is satisfied by the  occupation measures of planar Brownian  motion and SLE curves under natural parametrizations.
	Along the way we obtain general {positive} moment bounds for Gaussian multiplicative chaos.
	Contrary to previous results,  this does not require any assumption on the underlying measure $\sigma$ such as scale invariance, and hence may be of independent interest.
\end{abstract}

\section{Introduction} 
\subsection{Motivation and statement of the problem}
\label{SS:motiv}

Let $\gamma \in (0,2)$.  Liouville quantum gravity (LQG) is the random planar geometry whose volume measure, is given (formally) by
$$\mu_h=e^{\gamma h}dz,$$ where $h$ is a variant of the two-dimensional  Gaussian free field  (GFF), {and $dz$ is the Lebesgue measure on $\C$.}
The GFF is not a pointwise defined function but is rather a random distribution, so that making sense of {$\mu_h$ in the above definition} requires some nontrivial work.

{
This measure, now known as Liouville measure, can be constructed via a regularization and normalization procedure known as \textbf{Gaussian multiplicative chaos}; see \eqref{e.mudef}. This fact goes back in some form to the pioneering work of Kahane \cite{kahane} and (for a smaller range of values of the parameter $\gamma$) to the earlier work of H{\o}egh-Krohn \cite{hoeghkrohn}. It was subsequently rediscovered in \cite{KPZ} (see also \cite{RhodesVargas}) which formulated a framework for the construction of Liouville quantum gravity, as arising in the work of Polyakov \cite{polyakovstrings} in string theory. By the so-called \textbf{DDK ansatz} this geometry is also predicted to describe the scaling limit of natural random planar maps models. See \cite{cf:Da,DistKa,Nak,KPZ,BerestyckiKPZnotes,mot-survey, DKRV} for  physics and mathematics background on the LQG and \cite{BerestyckiKPZnotes, DKRV, mot-survey} in particular for the conjectured relation to random planar maps.

The construction of the LQG area measure was later complemented by further LQG observables, associated with Schramm--Loewner evolution \cite{Zipper,mating} and planar  Brownian motion \cite{lbm,diffrag}. In both cases we have a conformally invariant process whose trajectory is independent from the GFF but whose parametrization respects the LQG geometry.
	In the latter case, the reparametrized Brownian motion is called \textbf{Liouville Brownian motion}.
	More recently, a random metric associated with LQG was constructed in a remarkable series of works culminating in \cite{DDDF,GM}. (In the case where the parameter $\gamma=\sqrt{8/3}$, this metric coincides with the one constructed earlier in \cite{LQGandTBMI,LQGandTBMII,LQGandTBMIII}
	via an indirect mechanism called quantum Loewner evolution). All these observables have natural interpretations from the conjectured link with random planar maps.}

{The LQG area measure is perhaps the simplest and most well studied  LQG observable at this point.
	Moreover, much of the theory of LQG is based on the implicit assumption that this observable captures all of the geometry of LQG.
	For example, the \textbf{coordinate change formula} dictating whether two LQG surfaces should be viewed as conformally equivalent,
	is prescribed by the requirement that the LQG area measure must transform covariantly.
	Moreover, in the mating of trees theory (as introduced in \cite{mating}) a pair of Brownian motions $(L, R)$ is used to describe the LQG  area measure decorated by an SLE curve.
	Based on this assumption,  the LQG field and hence the entire underlying surface is captured by $(L, R)$; see Corollary \ref{C:mating}.}

This suggests that the Liouville measure is an observable that is sufficiently rich to capture
all the
relevant geometry. The main goal of this paper is precisely to put the assumption above on solid ground, by demonstrating that the underlying Gaussian free field $h$ can be measurably recovered from its Liouville measure (Theorem \ref{Theorem:measure determines field}). While this is intuitive at some level, this is also far from obvious a priori, since it is well known that the Liouville measure concentrates on points that are in some sense exceptional for the GFF, {which has Hausdorff dimension less than 2}.

In fact, we consider the same question in a more general (and, as we will see, more challenging) yet also very natural setup.
Let $\cA$ be a closed set and let $\sigma$ be a measure supported on $\cA$. We suppose that {$X$ is a variant of two dimensional GFF and}  $\sigma$ (and hence also $\cA$) is either fixed deterministic or random but independent from  $X$.  Let $d$ be the dimension of $\sigma$ (the maximal value such that $\sigma$ has finite $d$-dimensional energy) and let $\gamma < \sqrt{2d}$. Then $X$ induces a \textbf{Gaussian multiplicative chaos} measure $\mu_X$ supported on $\cA$ with parameter $\gamma$, that is, informally,
$$
d\mu_X = e^{\gamma X} d \sigma.
$$
In the special case where $X$ is a GFF with Dirichlet boundary conditions on a domain $D$ and $\sigma$ is the Lebesgue measure (henceforth the \textbf{Lebesgue case}) we recover the above setup.
In Theorem~\ref{thm:gen} we will show that under mild assumptions on $\sigma$ and on $X$, $\mu_X$ entirely determines $X$ restricted to $\cA$: more precisely (since $X$ is not defined pointwise), $\mu_X$ determines the harmonic extension of $X$ off $\cA$ (which is well defined under our assumptions).

As an application of this general framework, we obtain the following corollaries:
\begin{itemize}
	
	\item[(1)] a single trajectory of Liouville Brownian motion, viewed as a parameterized curve up to any given time $t\in [0, \infty]$, determines the Gaussian free field entirely on its range up to time $t$. In particular, if $t = \infty$ and the underlying domain is the whole complex plane, since the infinite range of Liouville Brownian motion is dense, that trajectory determines all of the GFF.

	\item[(2)] the quantum length of an SLE curve (as defined e.g. in \cite{KPZ} or \cite{mating}) determines the GFF on its range.

\end{itemize}

Along the way, we will obtain some moments estimates (both positive and negative) for Gaussian multiplicative chaos on such fractals which we feel are of independent interest, see Theorems \ref{thm:positive} and \ref{C:negmom} respectively.

\subsection{The equivalence theorem {for the LQG area measure}}\label{sec:main-results}
Let $D\subset\C$ be a {bounded} domain  and $h$ be a Gaussian free field with zero boundary condition (zero boundary GFF) on $D$. See \eqref{dirichlet prod} for its definition.
We can use a regularization procedure to define an area measure on $D$:
\begin{equation} \label{e.mudef}
	\mu_h := \lim_{\eps \to 0} \eps^{\gamma^2/2} e^{\gamma h_\eps(z)}dz,
\end{equation}
where $dz$ is Lebesgue measure on $D$, $h_\eps(z)$ is the mean value of $h$ on the circle $\{w\in\C: |w-z|=\eps \}$.
The limit represents weak convergence in the space of measures on~$D$, which exists in probability and is independent of the choice of regularization \cite{BerestyckiGMC}, \cite{Shamov}. In particular, $h$ determines $\mu_h$ almost surely.
In fact the convergence is almost sure for the above circle averages if $\eps$ is restricted to powers of two \cite{KPZ} (more recently, the restriction to powers of two was removed in \cite{SheffieldWang}).

Note that if $h = h_0 + g$, where $h_0$ is a zero-boundary GFF on $D$, and $g$ is a possibly random continuous function on $D$,
{the  LQG area measure $\mu_h$ associated to $h$ can still be defined as the almost sure limit from~\eqref{e.mudef}, hence is determined by $h$.}

Our first main result shows that the converse measurability holds.

\begin{theorem}\label{Theorem:measure determines field}
	Let $ h=h_0+g $ where $h_0$ is a zero boundary $ \GFF  $ on a bounded domain $ D\subset \C $ and $g$ is a random continuous function on $D$.
	Denote by $ \mu_h $ its LQG area measure with parameter $ \gamma\in (0,2) $.
	Then $ h $ is determined by $ \mu_h $ almost surely. That is, $h$ is measurable with respect to the $\sigma$-algebra generated by $\{\mu_h(A): A \text{ open  in }D\}$.
\end{theorem}
If $g \in {H_0^1(D)}$ is deterministic, then we note that {by the Cameron--Martin theorem} $h + g$ is absolutely continuous with respect to $h$ so in this case the theorem follows trivially from the case $g \equiv 0$. But here we only assume that $g$ is continuous, so $g$ can be much rougher. Moreover $g$ may depend on $h$.
\begin{remark}\label{rmk:generality}
	Theorem \ref{Theorem:measure determines field} covers various types of GFFs (Dirichlet boundary conditions, Neumann boundary conditions, mixed boundary conditions, the whole plane GFF, etc.) via the domain Markov property. Likewise, by absolute continuity,  Theorem \ref{Theorem:measure determines field} also covers the quantum surfaces defined in \cite{mating} including quantum cones, wedges, spheres and disks. {By~\cite{decomposition}, for a large class of log-correlated field, restricted to small enough neighborhood. the field can be written in the form of $h_0+g$ as in Theorem~\ref{Theorem:measure determines field}. Therefore this extends the result to fields considered in~\cite{decomposition}.}
\end{remark}

\paragraph{Application to the mating of trees theorem.} We now briefly explain an application of Theorem \ref{Theorem:measure determines field} to
the mating-of-trees framework to LQG developed in \cite{mating}. In the main result of that paper (Theorem 9.1), the authors consider a space-filling variant of SLE$_\kappa'$, $\kappa' = 16/\gamma^2$, on top of a $\gamma$-quantum cone, where the curve $\eta'$ is parametrized by its quantum area (i.e., $\mu_h (\eta'([s,t])) = t - s$ for all $s\le t \in \R$). We refer to \cite{mating} and \cite{Zipper} for the notion of quantum cone, while the space-filling variant of SLE was introduced in \cite{IG4}. The main theorem of \cite{mating} is that the left and right boundary quantum length of the curve $\eta'([0,t])$, relative to time 0, evolve as a certain two dimensional Brownian motion $(L_t,R_t)_{t\in \R}$ whose correlation coefficient is given by $\cos(\pi \gamma^2/4)$. (In fact, this formula was only proved for $\gamma \in [\sqrt{2}, 2)$ in \cite{mating}, and the corresponding result for $\gamma \in [0, \sqrt{2})$ is addressed in \cite{covariancestory}.)
In \cite[Chapter 10]{mating}, the authors proved that this pair of Brownian motions in fact determines the Liouville measure on the $\gamma$-quantum cone as well as the space-filling SLE almost surely, up to rotations, and used Theorem~\ref{Theorem:measure determines field} of this paper to conclude that this in turn determines the free field $h$ (up to rotations).
{More precisely, we have the following.
	\begin{corollary}\label{C:mating}
		Modulo rotations, the  field $h$ defining the $\gamma$-quantum cone is almost surely measurable with respect to the Brownian motion $(L_t,R_t)_{t \in \R}$.
	\end{corollary}}
	Other applications of Theorem~\ref{Theorem:measure determines field} (all subsequent to the time that a first version of this paper was made available on arxiv) can be found e.g. in~\cite{LQGandTBMIII,ALS,Dubedat-Shen}.

	\subsection{Equivalence theorem in a more general setup} \label{subsec::moregeneral}
	
	Recall that a Borel measure on a domain $D$ is {\em locally finite} if every point has a neighborhood of finite measure (or equivalently, if every compact set has finite measure). We will be interested in random pairs $(\sigma, \cA)$, where $\sigma$ is any (possibly random) locally finite measure on $D$ and $\cA$ is the (closed) support of $\sigma$. For example, $\cA$ could be one of the random fractal sets that arise in $\SLE$ theory, and $\sigma$ could be a `natural' fractal measure associated to $\cA$. Let $h$ be an instance of the GFF on $D$ with some boundary conditions chosen independently from $(\sigma, \cA)$.
	
	Fix $d \in (0,2]$ and assume that $\sigma$ has
	finite $(d-\eps)$-dimensional energy for all $\eps\in (0,d)$, i.e.,
	\begin{equation}\label{dim}
\cE_{d- \eps}: =		\iint \frac1{| x-y|^{d- \eps}}\sigma(dx) \sigma(dy) < \infty, \,\,\,\,\, \textrm{for all} \,\,\, \eps>0.
	\end{equation}
	The reader may recall that, by Frostman's theorem, the Hausdorff dimension of a closed set $\cA$ is the largest value of $d$ for which there exists a non-trivial measure $\sigma$ on $\cA$ satisfying \eqref{dim}. In the discussion below, we will not require that $d$ is the dimension of $\cA$, or that $\sigma$ is in any sense an optimal measure on $\cA$. Once $\sigma$ is fixed, choosing a smaller $d$ than necessary for \eqref{dim} will in some sense be equivalent to choosing a smaller $\gamma$.
	
	Let $\gamma <  \sqrt{2d}$. By Kahane's theory of multiplicative chaos (as explained, e.g., in Theorem 1.1 in \cite{BerestyckiGMC}) there is a way to define a measure $\measure$ (which depends on $h$ and $\sigma$) that can be formally written as:
	\begin{equation}\label{eq:GMC}
		{\measure (dz)} = \exp (\gamma h(z) - \frac{\gamma^2}2 \E(h(z)^2)) \sigma(dz).
	\end{equation}
	We will not explain the details of this construction here, but we emphasize that the measure $\measure$ is non-trivial,
	in the sense that $\measure (\cA) \in (0, \infty)$ a.s. and that its support is {also} $\cA$.
	
	We can now formulate the question we have in mind:
	
	\textbf{Question:} To what extent does the measure $\measure$ determine the field $h$?
	
	Clearly $\measure$ can only determine the field $h$ `restricted to $\cA$' in some sense. The issue of whether the restriction of $h$ to a fractal subset $\cA$ makes sense is itself not obvious. But if $\cA$ is deterministic (or more generally a `local' set coupled with $h$) then there is a natural way to define the harmonic extension (to the complement of $\cA$) of the values of $h$ on $\cA$; see \cite{localsetpaper}. (If a Brownian motion hits $\cA$ with probability 0, then this extension is just the {\em a priori} expectation of $h$.)
	In other words, we apply the \textbf{domain Markov property} {of $h$} (Theorem 1.26 in \cite{BerestyckiKPZnotes}) to write
	\begin{equation}\label{markov}
		h = \tilde h + h^{\har}
	\end{equation}
	where $\tilde h, h^{\har}$ are independent, $\tilde h$ is a {zero-boundary GFF  on $D\setminus\cA$} and $h^{\har}$ is harmonic in {$D\setminus\cA$}.
	
	In this paper we give the following condition on $(\cA,\sigma)$ ensuring that the measure $\measure$ a.s.\ determines the harmonic extension $h^{\har}$ of $h$ off $\cA$.
	For concreteness, we assume that $\sigma$ is a {deterministic} finite measure, $D$ is the open unit disk $\D=\{z\in\C:|z|<1\}$, and $\cA\subset \D$.
	In application to random fractals where $(\sigma,\cA)$ and $h$ are independent, we will check that the condition holds almost surely.

	We say that $(\cA, \sigma)$ satisfies {\bf Property (P)} if and only if there exist positive  constants  $c,r$  depending  on $(\cA, \sigma)$ such that
	\begin{equation}\label{eq:property}
		\min_{x\in \cA}\sigma(B_\eps(x)) \ge   c{\eps^{r}} \qquad\; \forall  \eps\in (0,1),
	\end{equation}
	where $B_\eps(z)$ be the Euclidean ball of radius $\eps$ centered at $0$.
	
	We say that $\cA$ satisfies {\bf Property (Q)} if and only if there exist positive  constants  $c,q$  depending  on $\cA$ such that
	\begin{equation}\label{eq:propertyQ}
		{\omega(z, B_\eps(x)) \le \frac{c\eps^q}{|z-x|^q}}\qquad\; \forall  \eps\in (0,1) \textrm{ and }x\in \cA, \; z\in \D \setminus\cA
	\end{equation}
	where $\omega(z, dx)$ be the harmonic measure of $\partial \cA$ viewed from $z$, and $\dist(z,\cA)$ is the Euclidean distance from $z$ to $\cA$. We will comment on these properties in Section \ref{subsec:examples}; in particular see Lemma \ref{lem:P} and Lemma \ref{lem:Q} for some simple conditions guaranteeing that these properties hold together with concrete examples of interest.

	\begin{theorem}\label{thm:gen}
		{Let $\sigma$ be a (deterministic) finite measure with closed support $\cA\subset \D$ such that  Properties (P) and (Q) hold.
			Let $h$ be  a zero boundary GFF on $\D$.  Given  $d\in(0,2]$ satisfying~\eqref{dim} and $\gamma\in (0,\sqrt{2d})$,
			the measure $\measure$ from~\eqref{eq:GMC} almost surely determines the function $h^{\har}$ from~\eqref{markov}.}
	\end{theorem}
	\begin{remark}\label{rmk:general2}
		{By the same reduction  in the proof of Theorem~\ref{Theorem:measure determines field},
			we can extend Theorem~\ref{thm:gen}  to the cases where the field  is a GFF plus a continuous function and the domain is general.}
	\end{remark}

Along the way, we derive general {positive} moments bounds for the mass of Gaussian multiplicative $\measure (\D)$.
Earlier works on this question (see e.g. \cite{rhodes-vargas-review}) assume either that $\sigma$ is the Lebesgue measure or, more generally, that it has nice scaling properties. The result we present below does not make any assumption on $\sigma$ beyond having finite $d-\eps$ energy as in \eqref{dim} (in particular, properties (P) or (Q) are not assumed).

\begin{theorem}\label{thm:positive_intro}
	{Fix $0<\mathbf{d}\le d$ such that $\cE_{\mathbf{d}} = \iint |x-y|^{- \mathbf{d}} \sigma(dx) \sigma(dy) < \infty$ and $0 \le \gamma< \sqrt{2 \mathbf{d}}$.}
	For each $\alpha \in (1,{\frac{2\mathbf{d}}{\gamma^2}} \wedge 2)$,  there exists a constant $c=c(\bf d,\gamma, \alpha)$ such that
	\[
	\E[{\measure}(\D)^\alpha] \le c\sigma(\D)^{3- \alpha} (\cE_{{\bf d}})^{\alpha-1}<\infty.
	\]
\end{theorem}
{Theorem~\ref{thm:positive_intro} is proved as Theorem~\ref{thm:positive} in Section~\ref{sec:GMC}. Bounds on $\E[{\measure}(\D)^\alpha]$ for negative $\alpha$ are also in~\cite{negative-moment}. However, in contrast to Theorem~\ref{thm:positive_intro}, the estimate in~\cite{negative-moment}  is only made explicit in terms of $\sigma(\D)$ and $\cE_{\bf d}$ when $\gamma\in (0,\sqrt{\bf d})$, namely, the $L^2$ regime; see \cite[Corollary~3.2]{negative-moment}.
In Section~\ref{sec:GMC} we use Theorem~\ref{thm:positive_intro} to give a  a similar estimate on negative moments  for \emph{all} $\gamma\in (0,\sqrt{2\bf d})$; see Lemma~\ref{lem:t} and Theorem~\ref{C:negmom}.}

	\subsection{Examples}\label{subsec:examples}
	{Property (P) is a quantification of the statement that $\sigma$ is sufficiently spread  out on $\cA$. It is rather weak since $q$ is allowed to be arbitrarily large. For example, it is satisfied if by occupation measures of Holder continuous curves.
		\begin{lemma}\label{lem:P}
			Suppose $\eta:[0,1] \mapsto \D$ is a Holder continuous curve. Namely, there exists positive constants $C,\alpha$ such that $|\eta(s)-\eta(t)|\le C|s-t|^\alpha$ for $s,t\in [0,1]$.
			Let $\cA=\eta([0,1])$. Let $\sigma$ be the occupation measure of $\eta$, which means for each Borel set $B\subset \D$, $\sigma(B)$ measures the amount of time $\eta$ spends in $B$.
			Then $(\cA,\sigma)$ satisfies Property (P) with $q=\frac1\alpha$.
		\end{lemma}
		\begin{proof}
			For $x\in \cA$, suppose $\eta(t)=x$ for some $t\in [0,1]$.
			Then $\eta(t')\in B_\eps(x)$ if $C|t-t'|^\alpha<\eps$. Therefore $\sigma(B_\eps(x))> (\eps/C)^{\frac1\alpha}$.
		\end{proof}
		The following sufficient condition for Property (Q) is a direct consequence of Beurling's estimate {(see \cite{Lawler-book} Equation (3.17)).}
		\begin{lemma}\label{lem:Q}
			Property (Q)  is satisfied for  $\cA$ with $q=\frac12$ as long as  $\cA$ is connected.
		\end{lemma}}

			\paragraph{Liouville Brownian motion.}
			For concreteness, let $\cB = (\cB_t, t \le T_D)$ be a (standard) Brownian motion starting from 0, run until it leaves $D$,
			{where $D$ is a domain such that $0\in D$ and $D\cup\partial D\subset \D$.}
			Let $\sigma$ denote the corresponding occupation measure: as above,
			for a Borel set $B \subset D$, $\sigma(B) = \int_0^{T_D} 1_{\{ \cB_s \in B\}} ds$. Its  closed support  $\cA$ is the range $\{B_t, t \le T_D\}$ of $\cB$. It is well known that the dimension of planar Brownian motion is  almost surely 2 and so condition \eqref{dim}  is a.s.\ satisfied by $\sigma$ with $d=2$.
			Furthermore, by Lemmas~\ref{lem:P} and~\ref{lem:Q}, both Properties~(P) and (Q) are satisfied by $(\sigma,\cA)$ almost surely, so that Theorem~\ref{thm:gen} applies.
			
			For $\gamma\in(0,2)$, let $h$ be {a zero-boundary GFF on $\D$}, independent of $\cB$.
			Following  \cite{diffrag,lbm}, the following limit makes sense in probability and defines an increasing process $(\phi(t), 0 \le t \le T_D)$, the \textbf{quantum clock}:
			$$
			\phi(t): = \lim_{\eps \to 0} \eps^{\gamma^2/2}\int_0^t e^{\gamma h_\eps(\cB_s)} ds.
			$$
			By definition, Liouville Brownian motion (LBM) is the reparametrization of $\cB$ by its quantum clock: namely,
			$$
			Z_t := \cB(\phi^{-1}(t)); \ \  \textrm{for } t \le \tau_D := \phi(T_D).
			$$
			The quantum clock can be rewritten in terms of the Gaussian multiplicative chaos $\measure$ associated to $h$ and $\sigma$ in the sense of \eqref{eq:GMC} as follows.
			For any $t \le T_D$,
			$$\phi(t)  = \int_D 1_{\{x \in \cB[0,t] \}} d\measure(x),$$
			namely, $\phi(t)$ is the $\measure$-mass of $\cB[0,t]$. Note that $(\cB[0,t], t \le T_D)$ is measurable with respect to $(Z_t,  t \le \tau_D)$: indeed it suffices to reparametrize $Z$ by its quadratic variation.
			Hence $(\phi(t), 0 \le t \le T_D)$ is a measurable function of $(Z_t, t \le \tau_D).$ As a consequence, we see that $\measure$ is a measurable function of $(Z_t, t \le \tau_D).$ Applying Theorem~\ref{thm:gen}, we deduce  that  the harmonic extension of $h$ off  the range of $\cB_{[0,T_D]}$ is a.s.\ determined by $(Z_t, 0 \le \tau_D)$.
			
			\medskip
			
			{\bf SLE.} Let $\kappa <8$ and let $\eta$ denote a chordal SLE$_\kappa$ from 0 to $\infty$ in the upper half plane $\mathbb{H}$. It is shown in \cite{LawlerRezaei} that an $\eta$ a.s.\ has nontrivial $d$-dimensional Minkowski content with $d=(1+\frac{\kappa}{8})\wedge 2$,
			which defines the so-called \emph{natural measure} on the curve. It is proved in \cite{LawlerRezaei} that SLE is H\"older continuous under its natural parametrization, for any $\kappa \in (0,8)$. (The optimal H\"older exponent is achieved in \cite{Zhan-holder}.) By Lemmas~\ref{lem:P}~and~\ref{lem:Q}, Property~(P)  holds almost surely for the occupation measure $\sigma$ of SLE under this natural parametrization, and Property~(Q) holds for an SLE curve segment.
			
			For $\kappa\in (0,4)$, it was proved by \cite{Benoist} that if $h$ is an independent Neumann GFF (with an additional logarithmic singularity at the origin), then the corresponding GMC measure with parameter $\gamma/2$ with respect to the reference measure $\sigma$, $e^{\frac{\gamma}2 h}\sigma(dz)$ is identical to the quantum length measure of $\SLE_\kappa$, as defined in \cite{Zipper}, where $\gamma = \sqrt{\kappa}$.
			The intricacy of~\cite{Benoist} lies in that it does not assume the existence of Minkowski content from \cite{LawlerRezaei}. Instead, it shows that if one takes the conditional expectation of the quantum length given the SLE curve,  this defines a measure over the range of $\eta$ which can be shown to satisfy the set of axioms in~\cite{LawlerSheffield}, which characterize the natural parametrization uniquely: hence this conditional expectation is $\sigma$.
			As a consequence, by Theorem \ref{thm:gen}, the quantum length of $\eta$ determines the restriction of $h$ to $\eta$ (or, more precisely, the harmonic extension $h^{\har}$ off of $\eta$).

			In fact, for $\kappa\in(4,8)$, the counterpart of \cite{Benoist} still holds, with the quantum length replaced by the so-called quantum natural time considered in \cite{mating}, and the factor $\frac12$ in the GMC adjusted according to the KPZ relation.
			
			\medskip
			
			\begin{remark}We expect that Property (Q) holds for  Cantor-like random fractals such as cut points of Brownian motion, cut points and double points of $\SLE_{\kappa}$ with $\kappa\in(4,8)$. Moreover, Property (P) holds for natural measures on them defined via the Minkowski content.
			\end{remark}
			
			\subsubsection*{Notations}
			The following notations are used throughout this paper. For $r>0,z\in \C$, we let $B_r(z)\subset\C$ be the ball of radius $r$ centered at $z$ and $\D=B_1(0)$, $r\D=B_r(0)$.
			The symbol $\eps$ always represents a small enough positive number (e.g. $\eps\in (0,1/100)$).  Suppose $h$ is a variant of GFF, let $h_\eps(z)$ be $h$ averaged along $\partial B_\eps(z)$.
			We let $C$ be constants arising in the argument that can vary line by line.
			If a function $\psi$ is such that $\psi(\eps)=o(\eps^{-p})$ for all $p>0$, we call $\psi$ a \emph{sub-polynomial} function. The main example of sub-polynomial functions are powers of $|\log\eps|$.
			\subsubsection*{Outline of the paper}
			In Section \ref{sec:preliminaries}, we provide some background on the Gaussian free field and prove a useful preliminary estimate. In Section~\ref{sec:LQM-GFF}
			we give the proof of Theorem~\ref{Theorem:measure determines field}. Although it will be a consequence of Theorem~\ref{thm:gen}, the proof is fairly elementary and contains the high level ideas of the more general case. In Section~\ref{sec:gen}; we prove Theorem~\ref{thm:gen}.

			\bigskip
			\noindent{\bf Acknowledgments.}
			We thank the Isaac Newton Institute for Mathematical Sciences, Cambridge, for generous support and hospitality during the programme \emph{Random Geometry} where part of this project was undertaken. The first and second authors were partially supported by EPSRC grants EP/GO55068/1 and EP/I03372X/1, and the first author by an FWF grant on ``Scaling limits in random conformal geometry". The second author was partially supported by a grant and a sabbatical fellowship from the Simons Foundation.
			The second and third authors were partially supported by NSF Award DMS 1209044. The third author was
			supported by a Junior Fellow award from the Simons Foundation and NSF Grants DMS-1811092 and DMS-2027986.
			Part of this work was undertaken during various visits by the first author to MIT and Columbia, respectively. Their hospitality is gratefully acknowledged.

\section{Proof of Theorem~\ref{Theorem:measure determines field}: the Lebesgue case}\label{sec:LQM-GFF}

\subsection{Preliminaries}\label{sec:preliminaries}
We start by briefly recalling some background on the Gaussian free field (GFF) with zero (sometimes also known as Dirichlet)
boundary conditions, before stating a useful estimate which will be used repeatedly in this paper.

Let $D$ be a domain in $\C$ with harmonically nontrivial boundary (i.e. the harmonic measure of $\partial D$ is positive as seen from any point in $ D $).
We denote by $H_0(D)$ the Hilbert-space closure of  $C_0^\infty( D)$ (the space of compactly supported smooth functions in $D$), equipped with the Dirichlet inner product
\begin{equation}\label{dirichlet prod}
	(f,g)_\nabla = \frac{1}{2\pi} \int_D \nabla f(z)  \cdot \nabla g(z) \,  dz.
\end{equation}
A zero boundary \textit{Gaussian free field}  on $D$  is given by the formal sum
\begin{equation}\label{h}
	h = \sum_{n=1}^\infty \alpha_n f_n,\quad (\alpha_n)\quad \textrm{i.i.d}\quad  N(0,1)
\end{equation}
where $\{f_n\}$ is an orthonormal basis for $H_0(D)$. Although this expansion of $h$ does not converge in $H_0(D)$, it can be shown that the convergence holds almost
surely in the space of distributions. See \cite{SheffieldGFF, BerestyckiKPZnotes} for more details.

The \textbf{domain Markov property} of the GFF (see Theorem 1.27 in \cite{BerestyckiKPZnotes}) states that given a subdomain $U \subset D$ and a zero boundary GFF $h$ on $D$, we can write
\begin{equation}\label{Markov}
	h = h^{\har} + h^{\supp}
\end{equation}
where $h^{\har}, h^{\supp}$ are independent; $h^{\har}$ restricted to $U$ is almost surely a harmonic function in $U$; and $h^{\supp}$ is a zero boundary GFF on $U$ extended to be zero on $D \setminus U$. Recall also that $h$ is \textbf{conformally invariant} (i.e., its law is invariant under conformal transformations of the complex plane -- see Theorem 1.30 in \cite{BerestyckiKPZnotes}). In particular, $h$ is scale-invariant.

We record a lemma which will be used frequently; this gives a quantitative control on the fluctuations of the harmonic part in a ball of radius $\eps$ arising in the Markov decomposition \eqref{Markov} when $U \subset D$ is taken to be a ball of radius $r >4\eps$.

\begin{lemma}\label{Lemma:modulus of harmonic extension}
	Let $D$ be a simply connected domain, and let $h$ be a zero boundary GFF on $D$. Suppose that $z \in D$ and $r>0$ is such that $U = B_r(z) \subset D$. Let $h^{\har}$ denote the harmonic part in the Markov decomposition \eqref{Markov} associated to $U$. For $\eps <r/4$, let
	$$
	\Delta_\eps(z) =\max\limits_{x\in  B_{\eps}(z)}h^{\har}(x) -\min\limits_{x\in  B_{\eps}(z)}h^{\har}(x).
	$$
	Then $ \E[\Delta^2_\eps(z)] \le C(\eps/r)^2 $ where $C$ is a universal constant independent of $\eps,r,z,D$.
\end{lemma}
\begin{proof}
	By translation and scaling invariance of the GFF, it suffices to prove the lemma in the case $z=0,r=1$. We first control the gradient of $ h^{\har} $ in $ B_{1/2}(0) $. In the proof of \cite[Lemma 4.5]{KPZ}, the authors show that the minimum of $ h^{\har} $ in $ B_{1/2}(0) $ has super exponential tail which is independent of the domain $D$ containing the unit disk.
	(In fact, we point out that a simpler proof of that lemma can be obtained using the Borell--Tsirelson inequality for Gaussian processes).
	The same is true for the maximum of $ h^{\har} $.  In particular,  the second moment of $ \| h^{\har} \|_{\infty, B_{1/2}(0)  }    $ is bounded by a universal constant $C$. By a standard gradient estimate of harmonic functions (see e.g. Lemma 2.8 in \cite{BerestyckiNorris}), $ \| \nabla h^{\har} \|_{\infty,B_{\frac{1}{4}}(0)}\leq C \| h^{\har} \|_{\infty, B_{1/2}(0)  }  $. So for $ \eps <1/4 $, we have
	$$\E[\Delta^2_\eps(z)]\leq C\eps^2 \E[\| \nabla h^{\har} \|^2_{\infty,B_{1/4}(0)}] \le C\eps^2,$$
	as desired.
\end{proof}

We will prove Theorem~\ref{Theorem:measure determines field} by taking the logarithm of the volume of a ball of radius $\eps$ (rescaled by $1/\gamma$) and showing that this does not deviate much from the circle average of the field at scale $\eps$ (up to a deterministic term). In other words, let us define $ h^{\eps} $ by setting
\begin{equation}\label{eq:estimator}
	e^{\gamma h^{\eps}(z)}= \mu_h(B_z(\eps)),\quad \textrm{i.e.,}\quad  h^\eps(z)= \gamma^{-1}\log\mu_h(B_z(\eps)).
\end{equation}
Roughly speaking, we will show that $ h^\eps -\E[h^\eps]$ converges to $ h $ in probability as $\eps \to 0$. Since $ h^\eps $ is determined by $ \mu_h $, so is $h$. We will achieve this via the following two lemmas which will be proved in Sections \ref{Section:proof of lemmas} and \ref{sec:cov}.

\begin{lemma}[Variance estimate]\label{Lemma:variance}
	Suppose $ D $ is a simply connected domain and $ D'=\{ x\in D| \dist(x,\partial D)>\eps_0  \}  $, where $ \eps_0 $ is a fixed constant. Suppose $ h $ is the zero boundary GFF on $ D $ and $ h^\eps $ is as in \eqref{eq:estimator}.  For all $z\in D', 0<\eps<\frac{\eps_0}{4}  $, let $ f_\eps(z)=h^\eps(z)-h_\eps(z) $. Then  we have
	\[
	\Var[f_\eps(z)]\leq C\log(\eps_0/\eps),
	\]
	where $ C $ is a universal constant independent of $D, \eps_0$ and $z$.
\end{lemma}

\begin{lemma}[Covariance estimate]\label{Lemma:covariance}
	Let h be the zero boundary $ \GFF  $ on $\D$ and $f_\eps $ be defined as in Lemma \ref{Lemma:variance} for $D=\D$. Then for $ x_1,x_2\in r\D $ and $|x_1-x_2|>\eps^{1/2}$,
	$$\Cov[f_\eps(x_1),f_\eps(x_2)]\le C_r\frac{\eps}{|x_2-x_1|}\log^{1/2}\frac{|x_1-x_2|}{\eps}$$where $C_r$ only depends on $r\in (0,1)$.
\end{lemma}

Given Lemma \ref{Lemma:variance} and Lemma \ref{Lemma:covariance}, we can get Theorem \ref{Theorem:measure determines field} in the case where $h$ is the zero boundary GFF on $\D$.
\begin{proposition}\label{Proposition: theorem for disk}
	If $ h $ is a zero boundary $ \GFF $ on $ \D $, then $ \mu_h $ determines  $ h $ almost surely.
\end{proposition}
\begin{proof}
	Suppose $\rho$ is a smooth function supported on $ r\D$ for some $r<1$. It is sufficient to show that $(h,\rho)$ is measurable with respect to $\mu_h$. We start by noting that
	\begin{align*}
		&\Var[(f_\eps,\rho)]=\int_{\D\times\D}dxdy \Cov[f_\eps(x),f_\eps(y) ]\rho(x)\rho(y)\\
		&=\int_{\{|x-y|>\eps^{1/2},x,y\in r\D\}}dxdy \Cov[f_\eps(x),f_\eps(y) ]\rho(x)\rho(y)+\int_{\{|x-y|<\eps^{1/2}\}}dxdy\Cov[f_\eps(x),f_\eps(y) ]\rho(x)\rho(y).
	\end{align*}
	
	By Lemma \ref{Lemma:covariance}$$1_{\{|x-y|>\eps^{1/2},x,y\in r\D\}}\Cov[f_\eps(x),f_\eps(y)]\le C_r\eps^{1/2}\log^{1/2}(\eps^{-1}).$$ Therefore
	$$\lim\limits_{\eps\to 0}\int_{\D\times\D}dxdy 1_{\{|x-y|>\eps^{1/2},x,y\in r\D\}}\Cov[f_\eps(x),f_\eps(y) ]\rho(x)\rho(y)=0.$$
	On the other hand, by Cauchy-Schwarz inequality and Lemma~\ref{Lemma:variance},
	$$\int_{\D\times\D}dxdy1_{\{|x-y|<\eps^{1/2}\}}\Cov[f_\eps(x),f_\eps(y) ]\rho(x)\rho(y)\le C\eps\log\frac{r}{\eps}.$$
	Therefore $\lim\limits_{\eps\to 0}\Var[(f_\eps,\rho)] =0$. Note that $\E((h_\eps, \rho)) = 0$. Hence
	\begin{equation}\label{L2pre}
		\lim_{\eps\to0}\;(h^{\eps},\rho)-\E[(h^{\eps},\rho)]-  (h_{\eps},\rho) = 0 \quad \textrm{in}\; L^2.
	\end{equation}
	
	To conclude, it remains to recall the standard fact that $h_\eps$ approximates $h$ in $L^2$ (when tested against $\rho$).
	\begin{lemma}\label{L:circleapprox}
		As $\eps \to 0$,  $ (h_\eps-h,\rho)$ converges to 0 in $ L^2(\P) $.
	\end{lemma}
	
	\begin{proof}
		Since $(h_\eps - h, \rho)$ is Gaussian, the result follows at once from the fact that $\cov (h_\eps(x), h_\eps(y))$ converges to $G(x, y)$ uniformly over compact subsets of $\D \times \D $ away from the diagonal, and from the bound
		\begin{equation}\label{eq:cov-circle}
			\cov (h_\eps(x), h_\eps(y)) \le - \log (|x-y | \vee \eps) + O(1)
		\end{equation}
		(see, e.g., Lemma 3.5 of \cite{BerestyckiGMC}).
	\end{proof}
	
	Combining \eqref{L2pre} and  Lemma \ref{L:circleapprox} we get that $(h^{\eps}  ,  \rho)-\E[(h^{\eps},\rho)]$ tends to $(h,\rho)$  in $ L^2 $.
	This implies that the random variable $(h,\rho)$ is measurable with respect to $ \mu_h $.
	
	So far, we have proved that for any smooth function $\rho$ supported on $\D$, $(h,\rho)$ is measurable with respect to $\mu_h$.  This yields that $h$ is determined by $\mu_h$ almost surely.
\end{proof}

With this in hand it is not hard to get a proof of Theorem~\ref{Theorem:measure determines field}.
\begin{proof}[Proof of Theorem \ref{Theorem:measure determines field}]
	We first assume $D=\D$ and write $h$ as $h_0+g$ where $ h_0 $ is an instance of a zero boundary $ \GFF  $ on $\D$ and $g$ is a random continuous function as in Theorem \ref{Theorem:measure determines field}.  Let $\mu_{h_0}$ be the Liouville quantum measure of $h_0$. Defined $h^\eps$ and $h^\eps_0$ by
	\begin{align*}
		e^{\gamma h^{\eps}}=\mu_h(B_\eps(z))\quad \textrm{and}\quad e^{\gamma h_0^{\eps}}=\mu_{h_0}(B_\eps(z)).
	\end{align*}
	
	By the intermediate value theorem, for all $x$, there is a $\xi_x$ such that
	\begin{equation}\label{eq:MVT}
		\mu_h(B_\eps(z))=\exp\{\gamma g(\xi_x)\}\mu_{h_0}(B_\eps(z)),\quad |\xi_x-x|\le\eps.
	\end{equation}
	Let $g^\eps(x)=g(\xi_x)$. By taking the logarithm on either side of \eqref{eq:MVT}, we have $h^{\eps}=h^{\eps}_0 + g^\eps$. Let $h_\eps, h_{0,\eps},g_\eps$ be the circle average process of $h,h_0,g$ respectively. Then $\forall \rho\in C^\infty_c(D)$,
	\begin{align}\label{eq:error}
		&(h^{\eps}  ,  \rho)-(h_\eps,\rho)=
		(h_0^{\eps} ,  \rho)-(h_{0,\eps},\rho)+ (g^\eps-g_\eps ,\rho).
	\end{align}
	By the argument in Proposition \ref{Proposition: theorem for disk}, $(h_0^{\eps} ,  \rho)-(h_{0,\eps},\rho)-\E[(h_0^{\eps} ,  \rho)]  $ tends to 0 in $L^2$ as $\eps$ tends to 0. Let $\omega_g$ be the modulus of continuity of $g$:
	\begin{equation}\label{eq:modulus}
		\omega_g(x,\eps)= \max\{|g(x)-g(y)|: |y-x|\le \eps \}, \quad \forall x\in D, \eps<\dist(x,\partial D).
	\end{equation}
	Since $g$ is a  continuous function,
	\begin{equation}\label{eq:def-loc}
		\lim_{\eps \to 0} (\omega_g(x,\eps),\rho(x) )=0\; a.s.\ , \quad \forall \rho\in C^\infty_c(D).
	\end{equation}

	By \eqref{eq:def-loc}, $(g^\eps-g_\eps ,\rho)$ tends to 0 a.s.\ as $\eps$ tends to 0. Hence
	\[
	\lim_{\eps \to 0} (h^{\eps} ,  \rho)-\E[(h_0^{\eps} ,  \rho)] =(h,\rho) \; \textrm{in probability}.
	\]
	As the same argument at the end of the proof of Proposition \ref{Proposition: theorem for disk}, we conclude the proof of Theorem \ref{Theorem:measure determines field} for $D=\D$.
	{By translation and scaling we get Theorem~\ref{Theorem:measure determines field} if $D$ is a ball in $\C$.}
	
	{For a general domain $D$, by the domain Markov property the field $h$ restricted to any ball $B\subset D$ can be written as a zero boundary GFF on $B$ plus a continuous function. Therefore  $h|_B$ is determined by $\mu_h$ almost surely. Varying $B$ we obtain
		Theorem \ref{Theorem:measure determines field}.}
\end{proof}

\subsection{The variance estimate}\label{Section:proof of lemmas}\label{sec:proof-Lemma}

\begin{proof}[Proof of Lemma \ref{Lemma:variance}]
	
	Consider first the disk $ U = B_{\frac{\eps_0}{2}} (z) $, and apply the Markov decomposition \eqref{Markov} to write $ h=h_0^{\supp}+h_0^{\har} $ as in that statement, where $ h_0^{\supp}  $ is a Dirichlet GFF on $U$  and $ h_0^{\har} $, restricted to $U$, is harmonic in $ U $ (and coincides with $ h $ on $ D \setminus U $).
	
	Let $ \xi_0 $ be a point in $B_{\eps} (z) $ such that
	$$
	\mu_{h}(B_{\eps} (z))= e^{\gamma h_0^{\har}(\xi_0) }   \mu_{h_0^{\supp}}(B_{\eps} (z)) .
	$$
	Since $h^{\har}_0$ has the mean value property, we have
	\begin{equation}\label{decomp1}
		f_\eps(z)=[ h_0^{\har}(\xi_0)-h_0^{\har}(z) ]+[\gamma^{-1} \log \mu_{h_0^{\supp}}(B_{\eps} (z))-h_{0,\eps}^{\supp}(z) ].
	\end{equation}
	By Lemma \ref{Lemma:modulus of harmonic extension},
	\begin{align*}
		\var[ h_0^{\har}(\xi_0)-h_0^{\har}(z)]  & \le \E [  (  h_0^{\har}(\xi_0)-h_0^{\har}(z) )^2 ] \\
		& \le \E ( \Delta_\eps(z)^2) \\
		& \le C(\eps/\eps_0)^2
	\end{align*}
	Applying the scaling and translation $ x\mapsto 2\cdot \frac{x-z}{\eps_0} $ to $h_0^{\supp}$ in \eqref{decomp1} (recall the conformal invariance of the Dirichlet GFF), in order to prove Lemma \ref{Lemma:variance} it suffices (by Cauchy--Schwarz) to show
	\begin{equation}\label{Ivar}
		\var(I) \le C |\log\eps|.
	\end{equation}
	where  \[ I= \tfrac1{\gamma}\log\mu_{h}(B_{\eps} (0))-h_{\eps}(0) \] and $ h $ is the zero boundary  $ \GFF  $ on $ \D $.

	We now start proving \eqref{Ivar}. We assume without loss of generality $ \eps=2^{-n} $ for some integer $n$, and recursively decompose $I$ as follows.
	Let  $ h^0=h$.
	For $0\le k<n-1  $, let $U_k = B_{2^{-k}} (0)$, and apply inductively the domain Markov property onto each $U_k$.
	That is, having defined $h^1, \ldots, h^k$ as Dirichlet GFFs on $U_1, \ldots U_k$ respectively, and $h^{\har}_1, \ldots, h^{\har}_k$ whose respective restrictions to $U_1, \ldots U_k$ are harmonic functions in the relevant domains, we apply the domain Markov property to $h^k$ and the subdomain $U_{k+1}$ to get a Dirichlet GFF $h^{k+1}$ on $U_{k+1}$ and an independent field $h^{\har}_{k+1}$ which is harmonic in $U_{k+1}$ such that
	$$
	h^k = h^{k+1} + h^{\har}_{k+1}.
	$$
	In this way, for each $k \ge 1$, $h^{k+1}$ and $h^{\har}_{k+1}$ are independent and more generally, $h^{k+1}$ is independent from the collection of random variables $h^{\har}_1, \ldots, h^{\har}_{k+1}$. Furthermore,
	$$
	h = \sum_{i=1}^{n-1} h^{\har}_i +  h^{n-1}
	$$
	and so taking the circle average at 0, by the mean value property of harmonic functions, we have
	\begin{equation}\label{dec_ca}
		\quad h_\eps(0) = \sum_{i=1}^{n-1} h_{i}^{\har}(0) + h^{n-1}_\eps (0).
	\end{equation}
	Let $ \xi_{k+1} $ be a point in $B_{\eps} (z) $ chosen in a measurable way with respect to $h^k$ such that
	$$\mu_{h^k}(B_{\eps}(0))=e^{\gamma h^{\har }_{k+1}  (\xi_{k+1}) } \mu_{h^{k+1}}(B_{\eps}(0)),$$
	whence
	\begin{equation}
		\label{dec_vol}
		\tfrac1{\gamma} \log \mu_h( B_\eps(0)) = \sum_{i=1}^{n-1} h^{\har} (\xi_{i}) + \tfrac1{\gamma}\log \mu_{h^{n-1}} (B_\eps(0)).
	\end{equation}
	Combining \eqref{dec_ca} and \eqref{dec_vol}, we can write
	\begin{equation}\label{Idec}
		I=\displaystyle{\sum_{k=1}^{n-1}} \Delta_k + R,
	\end{equation}
	where
	$$ \Delta_k=h_{k}^{\har} (\xi_k)-h_k^{\har}(0) \quad\textrm{and} \quad R=\tfrac1{\gamma} \log\mu_{h^{n-1}}(B_{\eps} (0) )-h^{n-1}_{\eps}(0) .$$
	The random variables $R$ and $\sum_{k=1}^{n-1}\Delta_k$ are not independent in general (because $\xi_k$ depends a priori on all of $h_k$). Nevertheless, by Cauchy--Schwarz,
	\begin{equation}\label{Idec2}
		\var(I)\le 2 \var (\displaystyle{\sum_{k=1}^{n-1}} \Delta_k)  + 2 \var (R).
	\end{equation}

	Mapping $B_{2^{-n+1}}(0)= B_{2\eps}(0) $ to $ \D $ by scaling and applying the LQG change of coordinates, we have
	\begin{equation}\label{eq:R}R\overset{\text{d}}{=} \tfrac1{\gamma} \log \mu_{h}(B_{\frac{1}{2}} (0) )-h_{\frac{1}{2}}(0) -\tfrac1{\gamma}Q\log 2\eps, \end{equation} where
	$ h $ is the zero boundary $ \GFF  $ on $ \D $, and $h_{\frac{1}{2} }(0)$ refers to the circle average of $\partial B_{1/2}(0)$. Since $\mu_h(B_{1/2}(0))$ has moments of  positive and negative orders (see Theorems 2.11 and 2.12 in \cite{rhodes-vargas-review} for the case $\gamma <2$, and Corollary 6 of \cite{DRSV1} for the case $\gamma =2$), we deduce that its logarithmic moments are also finite. Therefore $\Var[R]$ is a constant independent of $\eps$.

	Furthermore, by Lemma  \ref{Lemma:modulus of harmonic extension} and the definition of $ \Delta_k $, $ \E[\Delta^2_k]\leq C 4^{k-n} $ where $ C $ is independent of $ \eps $. Thus by the Cauchy--Schwarz inequality,
	$$\var (\displaystyle{\sum_{k=1}^{n-1}} \Delta_k)  \le \E[(\sum\limits_{k=1}^{n-1}\Delta_k )^2]\le n\E[\sum\limits_{k=1}^{n-1}\Delta^2_k]\le C\log\eps^{-1}$$
	
	Together with \eqref{Idec2}, this proves \eqref{Ivar} and hence Lemma~\ref{Lemma:variance}.
\end{proof}

\subsection{The covariance estimate}\label{sec:cov}

\begin{proof}[Proof of Lemma \ref{Lemma:covariance}]
	For fixed $ x_1 $ and $ x_2 $ in $ r\D $, let $ L $ be the line segment orthogonally bisecting $\overline{x_1x_2}$.
	Let $ U_i (i=1,2)$ be the connected component of $ \D\setminus L$
	containing $x_i$. Apply the domain Markov property in $U_1$ and $U_2$ to write
	$$ h = h^{\har}+ h^1+ h^2 $$
	where $h^{\har}$ is harmonic on $U_1$ and $U_2$, $h^1$ and $h^2$ are GFF (with Dirichlet boundary conditions) in $U_1$ and $U_2$ respectively, and all three terms are mutually independent. Then for $|x_1 - x_2|\ge \eps^{1/2}>2\eps$, the mean value property of $h^{\har}$ gives  	
	\[f_\eps(x_i) =\gamma^{-1} \log \int_{B_\eps(x_i)}e^{ \gamma h^{\har}(z) }  d\mu_{h^i} (z) -h^{\har}(x_i)-h_\eps^i(x_i).\]
	
	As in the proof of Lemma \ref{Lemma:variance}, let $ \xi_i $ be a point in $ B_\eps(x_i) $ (chosen in a measurable way with respect to $h^i$ and $h^{\har}$) such that $$\int_{B_\eps(x_i)} e^{ \gamma h^{\har}(z) }  d\mu_{h^i}=e^{ \gamma h^{\har}(\xi_i) } d\mu_{h^i}B_\eps(x_i).$$
	
	Let
	\[
	\Gamma_\eps(x_i)=\gamma^{-1} \log\mu_{h^i}B_\eps(x_i)- h_\eps^i(x_i)\quad \textrm{and}\quad \Delta_i=h^{\har}(\xi_i)-h^{\har}(x_i).
	\]
	Then $ f_\eps(x_i) =  \Gamma_\eps(x_i) +\Delta_i  $. Therefore
	$$ \Cov[f_\eps(x_1) , f_\eps(x_2)]= \Cov[\Gamma_\eps(x_1) , \Gamma_\eps(x_2)]+ \Cov[\Gamma_\eps(x_1) , \Delta_2]+ \Cov[\Gamma_\eps(x_2) , \Delta_1]+ \Cov[\Delta_1, \Delta_2].$$
	
	By Lemma \ref{Lemma:variance}, $\Var[\Gamma_\eps(x_i)]\le C_r\log\frac{|x_1-x_2|}{\eps}$. By Lemma \ref{Lemma:modulus of harmonic extension}, $ \Var[\Delta_i]\le C\frac{\eps^2}{|x_2-x_1|^2}$. $ \Gamma_\eps(x_1) $ and $ \Gamma_\eps(x_2) $ are independent, which means $\Cov[\Gamma_\eps(x_1) , \Gamma_\eps(x_2)]=0$. Therefore, by Cauchy--Schwarz inequality we have,
	\begin{align*}
		\Cov[f_\eps(x_1) , f_\eps(x_2)]&\le (\Var[\Gamma_\eps(x_1)]\Var[\Delta_2])^{1/2}+(\Var[\Gamma_\eps(x_2)]\Var[\Delta_1])^{1/2}+(\Var[\Delta_1]\Var[\Delta_2])^{1/2}\\
		&\le C_r\frac{\eps}{|x_2-x_1|}\log^{1/2}\frac{|x_1-x_2|}{\eps}.
	\end{align*}
	This concludes the proof of Lemma \ref{Lemma:covariance}, and with it the proof of Theorem \ref{Theorem:measure determines field}.
\end{proof}

\section{Moments on GMC over fractals}\label{sec:GMC}

As should be clear from the previous section, getting bounds (both positive and negative moments) for Gaussian multiplicative chaos with respect to the reference measure $\sigma$ supported on the (typically fractal) set $\cA$ is a key part of the argument. Such moments are usually derived from scale invariance consideration and Kahane's inequality. Clearly, such arguments are not directly applicable when the measure $\sigma$ is not itself assumed to be scale-invariant in any reasonable sense. We hence need to develop arguments to control these moments. We believe these results are of independent interest.

In this section, we let $h$ be the zero boundary GFF on {$D$}. We suppose that $D$ is bounded; after Lemma \ref{lem:LDP} we will take $D = \D$.
Let $\sigma$ satisfy \eqref{dim} {and $\sigma(D) < \infty$}. Let $0 \le \gamma < \sqrt{2d}$ (so we only consider the \textbf{subcritical} case). Fix a {${\bf d} \in (\gamma^2/2,d)$ so that}
such that $\gamma < \sqrt{2 {\bf d}}$ and
\begin{equation}\label{eq:energy}
	\cE_{{\bf d}} :=\iint_{D^2} |x-y|^{-{\bf d}} \sigma(dx) \sigma(dy) < \infty.
\end{equation}
Let $\mu_h=\measure$ be defined as in \eqref{eq:GMC} for this parameter $\gamma$.
\begin{theorem}\label{thm:positive}
	For each $\alpha \in (1,{\frac{2\mathbf{d}}{\gamma^2}} \wedge 2)$,  there exists a constant $c=c(\bf d,\gamma, \alpha)$ such that
	\[
	\E[\mu_{h}(D)^\alpha] \le c\sigma(D)^{{3- \alpha}} (\cE_{{\bf d}})^{\alpha-1}<\infty.
	\]
\end{theorem}

\begin{proof}
	Let $\bar x , \bar y$ be two  i.i.d. samples from $\sigma(\cdot \cap D) / \sigma(D)$. Note that $\cE_{{\bf d}}={\sigma(D)^{{2}}} \E ( | \bar x - \bar y |^{- {\bf d}} )$.
	
	Set $\delta=\alpha-1 \in (0,1)$, and write
	$$
	\E [ \mu_h({D})^{\alpha} ] = \E[\mu_h(D) \mu_h(D)^\delta] = {\E[\mu_h(D)]}\E^* [ \mu_h(D)^\delta]={\sigma(D)\E^*[\mu_h(D)^\delta]}
	$$
	where $\P^*$ denote the law of the field biased by $\mu_h(D)$.  Using Girsanov's theorem (see Lemma 2.5 in \cite{BerestyckiKPZnotes}), we can rewrite this as
	$$
	\E^*[\mu_h(D)^\delta] = \int_D \sigma(dx) \E( \left(\int_D e^{ \gamma G(x,y) }  \mu_h (dy) \right)^\delta)
	$$
	where $G(x,y) = G_{D} (x,y)$ is the Green function in $D$ with Dirichlet {(zero)} boundary conditions, normalised so that
	\begin{equation}\label{green}
		G(x, y) = - \log |x- y | + O(1) \text{ as }y \to x.
	\end{equation}
	For any integer $n\ge 0$, let $A_n(x)$ denote the annulus at distance $2^{-n}$ from $x$, i.e., $A_n(x) = \{ y: |y-x| \in [2^{-n-1}, 2^{-n})\}$. Then using \eqref{green} and the fact (coming from concavity of the function $t >0 \mapsto t^\delta$) that $(a_1 + \ldots + a_n)^\delta\le a_1^\delta + \ldots + a_n^\delta$ for any $0< \delta <1$ (i.e., $\alpha<2$ as per our assumptions) and $a_i >0$,
	we see that for some constant $C$ depending on $\delta$ and $\gamma$,
	\begin{align*}
		\E^*[\mu_h(D)^\delta]  & \le  C\sum_{n} \int_D \sigma(dx)  \E (  \left(\int_{A_n(x)} e^{-\gamma^2 \log |x-y| } \mu_h(dy) \right)^\delta) \\
		& \le C \sum_n 2^{n \gamma^2\delta}  \int_D \sigma(dx) \E ( \mu_h(A_n(x))^\delta).
	\end{align*}
	For each fixed $x$, consider instead of $h$ an exactly scale invariant field $X$ centered around $X$, i.e., so that
	$$\{X(x+\lambda z)\}_{z \in B_1(0)}=\{\tilde X(z)+\Omega_{-\log \lambda}\}_{z \in B_1(0)},
	$$
	where $\Omega_{r}$ is a Gaussian with variance $r$ independent of $\tilde X$.
	Set $\lambda = 2^{-n}{\le} 1$. Note that for $\delta <1$, and fixed $n$, making the change of variables $y = x + \lambda z$, and denoting $\sigma_{\lambda, x} (dz)$ the corresponding image measure (so that $y \in A_n(x)$ means $z \in A_1(0)$ and the total mass is $\sigma_{\lambda, x} (A_1(0)) = \sigma(A_n(x))$).
	\begin{align*}
		\E \left[ \left( \int_{A_n(x)} e^{\gamma X_{\lambda \eps} (y)} (\lambda \eps)^{\gamma^2/2} \sigma(dy) \right)^\delta \right]&\le  \lambda^{\gamma^2\delta /2} \E [ ( \int_{{A_1(0)}} e^{\gamma X_{ {\eps}} ({x+}\lambda z)} \eps^{\gamma^2/2} \sigma_{\lambda, x} (dz) )^\delta]\\
		& =   \lambda^{\gamma^2\delta /2} \E [ e^{\delta \gamma \Omega_{\log 1/\lambda}} ( \int_{A_1(0)} e^{\gamma \tilde X_\eps({x+}z)} \eps^{\gamma^2/2} \sigma_{\lambda, x} (dz))^\delta ]\\
		& = \lambda^{\gamma^2 \delta /2} e^{\delta^2 \gamma^2 \log (1/\lambda)/2} \E [  ( \int_{A_1(0)} e^{\gamma X_\eps({x+}z)} \eps^{\gamma^2/2} \sigma_{\lambda, x} (dz))^\delta ]\\
		& \le \lambda^{\gamma^2\delta/2 - \delta^2 \gamma^2 /2 }  \E [  ( \int_{A_1(0)} e^{\gamma X_\eps({x+}z)} \eps^{\gamma^2/2} \sigma_{\lambda, x} (dz)) ]^\delta,
	\end{align*}where the last inequality is by Jensen's inequality since $\delta \le 1$. As a consequence, by Kahane's inequality (Theorem 3.12 in \cite{BerestyckiKPZnotes}), there exists a constant $C>0$ depending only on $\gamma$ and $\delta$ but not $n$ or $\eps$, such that
	\begin{align*}
		\E [\mu_{\eps 2^{-n} } (A_n(x))^\delta] & \le C \lambda^{\gamma^2\delta/2 - \delta^2 \gamma^2 /2 }  \sigma(A_n(x))^\delta.
	\end{align*}
	Letting $\eps\to 0$, we get that for any ${n\ge 0}$,
	$$
	\E[\mu(A_n(x))^\delta] \le C 2^{- n (\gamma^2\delta/2 - \delta^2 \gamma^2 /2 )} \sigma(A_n(x))^\delta.
	$$
	We deduce
	$$
	\E^*[\mu(D)^\delta] \le C \sum_n 2^{n \gamma^2 \delta -  n (\gamma^2\delta/2 - \delta^2 \gamma^2 /2 ) } \int \sigma(dx) \sigma(A_n(x))^\delta.
	$$
	Now observe that
	\begin{align*}
		\sigma(A_n(x))&  \le \sigma(D)\P [|\bar y- x |  \le 2^{-n}]
	\end{align*}
	so that by Jensen's inequality again (since $\delta \le 1$)
	\begin{align*}
		\int_A\sigma(dx) \sigma(A_n(x))^\delta & \le  \sigma(D)^{{1+\delta}} \E ( \P( |\bar x-\bar y| <  2^{-n} | {\bar x } )^\delta)\\
		& \le \sigma(D)^{{1+\delta}}\P( |\bar x-\bar y| \le 2^{-n} )^\delta\\
		& \le  \sigma(D)^{{1+\delta}}\frac{ \E ( |\bar x-\bar y|^{-\bf d})^\delta}{2^{n{\bf d} \delta}},
	\end{align*}
	where in the last inequality we have used Markov's inequality.
	By choice of $\mathbf{d}$, the numerator in right hand side is finite.
	Now $\delta < 2{\bf d}/\gamma^2 -1$ implies that $\delta(  \gamma^2/2 - {\bf d})  + \gamma^2\delta^2/2 <0$.
	Putting everything together, we can find $c=c(d,\gamma,\delta)$ such that
	$$
	\E^*[\mu(D)^\delta] \le c \sigma(D)^{{1+\delta}} \E[|\bar x-\bar y|^{-\bf d}]^\delta.
	$$
	{Since $\cE_{{\bf d}}=\sigma(D)^{{2}} \E ( | \bar x - \bar y |^{- {\bf d}} )$ and $\delta=\alpha-1$,} we conclude the proof.
\end{proof}

We now turn to negative moments, and take $D = \D$. The key input will be the following statement taken from {\cite[Lemma~3.1]{negative-moment}}.
\begin{lemma}\label{lem:LDP}
	Set $\bar \beta=\max\{{\bf d},\sqrt{2{\bf d}}\gamma \}$.  Fix any choice of $\delta>0$,
	let
	$$
	{\beta=\frac{1+\delta}{1+2\delta} \bar\beta+ \frac{\delta}{1+2\delta}\gamma^2\in (\gamma^2,\bar \beta)\quad  \textrm{and}\quad \ell = \frac{\beta-\gamma^2}{ \beta+\gamma^2\delta}=\frac{\bar\beta-\gamma^2}{\bar \beta+2\gamma^2\delta}\in(0,1).}$$
	Let  ${\phi_\beta( x)}=\int_\D  |x-y|^{-\beta} \mu_h(dy)$ denote the (quantum) energy function with power exponent $\beta$ and let $\bar x$ be a point sampled proportionally to the measure $\sigma(dx)$.  Then {\(\P[\phi_{\beta}(\bar x)<\infty]=1\) and}
	\begin{align}\label{eq:probability}
		\E\left[\exp\left (-t\mu_h(\D)\right ) \right] \leq \frac{{32}}{\sigma(\D)t^\ell}  \qquad\textrm{if}\qquad 	\P\left (\phi_{\beta}(\bar x){\le} {2^{-4(1+\delta \ell)} t^{\delta\ell}} \right )\ge 1/2.	
	\end{align}
\end{lemma}	
We will apply the above lemma to get the following estimate about the tail of the Laplace transform of $\mu_h(\D)$, which controls the small ball probability of $\mu_h(\D)$ and hence its negative moments.
\begin{lemma}\label{lem:t}
	For $\delta$ and $\ell$ as  in Lemma~\ref{lem:LDP}, {and $\alpha = (2\mathbf{d}/\gamma^2)\wedge 2$}, there exists ${c_*=c_*}(\bf d,\gamma,\delta)$ such that
	\begin{equation*}
		\E\left[\exp\left (-t\mu_h(\D)\right ) \right] \leq \frac{{32}}{\sigma(\D)t^\ell}  \qquad\textrm{if}\qquad 	t^{\delta \ell} \ge c_*\cE_{\bf d}\sigma(\D)^{{-1 + \tfrac1{\alpha-1}}}.
	\end{equation*}
\end{lemma}

\begin{proof}
	Note that  $\beta/\gamma < \sqrt{2 {\bf d}}$ by our assumption on $\beta$. Let $\mu^{\gamma^{-1}\beta}_{h}$ be the Gaussian multiplicative chaos with respect to the reference measure $\sigma$ on the set $\cA$, with parameter $\gamma^{-1}\beta$.
		By Girsanov's theorem (see Lemma 2.5 in \cite{BerestyckiKPZnotes}),  for $r\in(0,1)$
	\begin{align*}
		\E[\phi_\beta(\bar x)^{r}] & \le c  \E [ \int_\D \left( \int_{\D}e^{\beta G(x,y)} \mu_h(dy) \right)^r \frac{\sigma(dx)}{\sigma(\D)} ]\\
		& \le c \sigma(\D)^{-1} \E\left[(\mu_h(\D))^{r}\mu^{\gamma^{-1}\beta}_{h}(\D)\right].
	\end{align*}
	{Since $\beta > \gamma^2$ and $\alpha = (2\mathbf{d}/\gamma^2) \wedge 2$, we have $\alpha \in (1,2{\bf d}\gamma^2\beta^{-2}\wedge 2)$.}
	Set $r=1-\alpha^{-1}$. By H\"older's inequality,
	\begin{align*}
		\E[\phi_\beta(\bar x)^{r}]   & \le \sigma(\D)^{-1}\E\left[\mu_h(\D)\right]^{r} \E\left[\left(\mu^{\gamma^{-1}\beta}_{h}(\D)\right)^{\alpha}\right]^{\alpha^{-1}}.
	\end{align*}
	Observe that
	{$\E\left[\mu_h(\D)\right]^{r} $ is bounded by $\sigma(\D)^r$ and $\E\left[\left(\mu^{\gamma^{-1}\beta}_{h}(\D)\right)^{\alpha}\right]$}
	can be bounded by Theorem~\ref{thm:positive}. Therefore  there exists ${c_0}=c_0({\bf d},\gamma,\delta,r)$ such that
	\[
	\E[\phi^{r}_\beta(\bar x)]\le c_0 \sigma(\D)^{{\frac{2}{\alpha}-1}} (\cE_{\bf d})^r=c_0 \sigma(\D)^{{\frac{1}{\alpha}-r}} (\cE_{\bf d})^r.
	\]
	Now let $t>0$ be such that {$t^{\delta \ell r}\ge 2^{1+4(1+\delta \ell)} \times c_0 \sigma(\D)^{{\frac{1}{\alpha}-r}} (\cE_{\bf d})^r$}. Then by Markov's inequality
	\[
	\P\left[\phi_{\beta}(\bar x)\ge {2^{-4(1+\delta \ell)} t^{\delta\ell}}  \right]\le {2^{4(1+\delta \ell)}\frac{\E[\phi^{r}_\beta (\bar x)]}{t^{\delta \ell r}}\le 1/2}.
	\] {Since  $\alpha r = \alpha -1$}, Lemma~\ref{lem:t} follows from Lemma~\ref{lem:LDP}.
\end{proof}

{As a consequence, we obtain negative moments for $\mu_h(\D)$.}
\begin{theorem}\label{C:negmom}
	{Set $\ell=\frac{\bar \beta-\gamma^2}{\bar  \beta+2\gamma^2}$ with $\bar \beta=\max\{{\bf d},\sqrt{2{\bf d}}\gamma \}$, {and suppose $\sigma(\D)\le 1$}.}
		For any $s\in (0,\ell)$ and $\cE' \ge \cE_{\bf d}$, there exists a constant $c=c(\cE',  s, \bf d, \gamma)$ such that
	\[
	\E[\mu_h(\D)^{-s}] \le {c(\cE',  s, \bf d, \gamma)} \sigma(\D)^{- s / \ell}<\infty.
	\]
\end{theorem}
{The assumption that $\sigma(\D) \le 1$ is to ensure the appearance of $\sigma(\D)^{- s / \ell}$ in the upper bound.
	By considering the restriction of $\sigma$ to a small enough ball, the statement $\E[\mu_h(\D)^{-s}]<\infty$ holds without this assumption.}
\begin{proof}[Proof of Theorem~\ref{C:negmom}]
	Note that $\int_{0}^{\infty} t^{s-1}e^{-t\mu_h(\D)}dt= \Gamma(s)\mu_h(\D)^{-s}$, where $\Gamma(\cdot)$ is the Gamma function.
	{Set $\delta=1$, and take $\ell=\frac{\bar \beta-\gamma^2}{\bar  \beta+2\delta\gamma^2} $ as in Lemma \ref{lem:LDP},} {and let $\cE' \ge \cE_{\bf d}$}. Since $\sigma(\D) \le 1$ and $\alpha = (2 \mathbf{d}/ \gamma^2) \wedge 2 >1$, so that $1/ (\alpha -1)>0$, we see that $\sigma(\D)^{-1 + 1/ (\alpha -1)} \le \sigma(\D)^{-1}$. Hence  the estimate of Lemma \ref{lem:t} is valid as soon as $t^\ell \ge c_*\cE' /\sigma(\D)$.
We therefore obtain the following bound:
	\begin{align*}
		\E[\mu_h(\D)^{-s}] & \lesssim \int_0^{({c_*} \cE' / \sigma(\D))^{1/\ell}} t^{s-1} dt + \int_{(c_* \cE' / \sigma(\D))^{1/\ell}}^\infty \frac{{32}}{\sigma(\D)} t^{s-1 - \ell} dt \\
		&
		\lesssim \frac1s\left(\frac{\cE'}{\sigma(\D)} \right)^{\frac{s}{\ell}} + \frac{1}{\sigma(\D)}\left(\frac{\cE' }{\sigma(\D)} \right)^{\frac {s- \ell}{\ell} }\\
		& \lesssim \sigma(\D)^{- s / \ell},
	\end{align*}where the implicit constants only depend on $\cE',  s, \bf d, \gamma$.
\end{proof}
\begin{remark}
	{Notice that this argument somewhat surprisingly relies on first establishing positive moments to obtain negative moments.
		In \cite{negative-moment}, the negative moments of \emph{all} order was proved for $\mu_h(\D))$. However, the estimate is only made explicit in terms of $\sigma(\D)$ and $
		\cE_{\bf d}$ when $\gamma\in (0,\sqrt{\bf d})$, namely, the $L^2$ regime. See \cite[Corollary~3.2]{negative-moment}. Our Lemma~\ref{lem:t} holds for all $\gamma\in(0,\sqrt{2\bf d})$.}
\end{remark}

\section{Proof of Theorem~\ref{thm:gen}}\label{sec:gen}

{Now suppose that $\cA\subset \D$ is a compact closed set. Since $\sigma$ is usually clear from the context, we simply write $\measure$ as $\mu_h$.}

For $z\in \D$, recall that $\omega(z,dx)$ is the harmonic measure of $\D\setminus \cA$ viewed from $z$.
For $\eps\in (0,1/2)$,  let $h_\eps$ be the circle average process of $h$ at $z$. Note that $h_\eps=0$ on $\partial \D$.
Let \begin{align}
	\tilde{h}^{\har}_\eps(z)&= \int {\omega}(z,dx) \left[\gamma^{-1} \log \mu_h( B_\eps(x)) \right],\nonumber\\
	h^{\har}_\eps(z)&= \int {\omega}(z,dx) h_\eps(x),\nonumber\\
	g_\eps(z)&= \int k_\eps(x) \omega(z,dx)=\tilde{h}^{\har}_\eps(z)-h^{\har}_\eps(z), \label{eq:k}
\end{align}where  $k_\eps(x)=0$ on $\partial \D$ and
\begin{align}\label{eq:def-k}
	k_\eps(x)=\gamma^{-1} \log \mu_h(B_\eps(x))-h_\eps(x),\quad \forall x\in \cA.
\end{align}
Let $h^{\har}$ be the harmonic extension of $h$ off $\cA$. Following the strategy in Section~\ref{sec:LQM-GFF}, we will show that $h^{\har}_\eps$ is a good estimator of $h^{\har}$ (Lemma~\ref{lem:harmonic-circle}) and $g_\eps$ tends to 0 (Lemma~\ref{lem:finite-variance} and~\ref{lem:finite-covariance}).
Before going into the proof, we first record a useful Brownian motion fact.

{In Theorem~\ref{thm:gen} we assume  that $\mathcal \cA\subset  \D$  satisfies Property (Q).
	This assumption is used to justify the intuitively obvious claim that $h^{\har}_\eps$ is a good estimator for $h^{\har}$.
	We will see from the following proof that a much weaker condition than Property (Q) would suffice.}

\begin{lemma}\label{lem:harmonic-circle}
	$\forall \rho\in C_0^\infty(D)$, $\E[(h^{\har}_\eps,\rho)]=0$ and $(h^{\har}_\eps,\rho)$ tends to $(h^{\har},\rho)$ in probability.
\end{lemma}
\begin{proof}
	Without loss of generality, we can assume $\rho$ is a continuous probability density function. A general $\rho$ can be write as $c_1\rho_1-c_2\rho_2$ where $\rho_1,\rho_2$ are such functions. Now we sample $z$ according  to $\rho(z)dz$ and then sample a Brownian motion from $z$. Suppose $X$ is the exit location of the domain $\D\setminus \cA$ for this Brownian motion and $U$ is a uniform unit 2D vector independent of $X$. Let $\hat \rho(x)dx$ and $ \hat\rho_\eps(x)dx$ be the distribution of $X$ and $X+\eps U$ respectively. Then $(h^{\har},\rho)=(h,\hat \rho)$ and
	$( h^{\har}_\eps,\rho)=(h,\hat\rho_\eps)$.  Since $h\overset{d}{=}-h$, we have $\E[(h^{\har}_\eps,\rho)]=0$.
	It suffices to prove that $\Var[(h,\hat \rho_\eps-\hat \rho)]\to 0$.
	
	Suppose $X$ and $Y$ are two independent copies sampled from $\hat{\rho}(z)dz$, and $U_1,U_2$ are two independent uniform unit vector. Then it suffices to show that
	\begin{equation}\label{Ggoal}
		\lim_{\eps\to 0}\E[G_{\D}(\bullet,*)]=\E[G_\D(X,Y)]
	\end{equation}
	where $\bullet$ represents either $X$ or $X+\eps U_1$, $*$ represents either $Y$ or $Y+\eps U_2$, and $G_{\D}$ is the Green function in $\D$. Here we extend $G_\D$ to $\C^2$ by setting $G_{\D}=0$ on $\C^2\setminus \D\times \D$.
	
	\eqref{Ggoal} has to be shown in all four possible combinations of the terms. We only explain the case $G_\D( X +\eps U_1, Y + \eps U_2)$, since  other cases are similar and simpler.
	
	By the mean value property of  $G_\D$  off the diagonal, we have
	$$
	\E[G_\D(X+ \eps U_1,Y + \eps U_2) \indic{|X-Y|> 2\eps,X,Y\in \cA}] = \E[G_\D (X, Y) \indic{|X-Y|> 2\eps,X,Y\in \cA}].
	$$
	On the other hand, by~\eqref{eq:cov-circle}, we have
	$$
	\E[ G_\D( X +\eps U_1, Y + \eps U_2) \indic{| X- Y | \le 2 \eps}]
	\le  C|\log\eps| \P[|X- Y | < 2 \eps].$$
	It suffices to show that $\lim_{\eps \to 0}|\log\eps| \P[|X- Y | < 2 \eps]=0$.
	Since $\cA$ satisfies Property~(Q), we have
	\[
	\P[|X- Y | < 2 \eps]\le \max_{x\in \cA} \P[|Y-x|<2\eps] \le \max_{x\in \cA}\int_{\D} \rho(z) \omega(z,B_{2\eps}(x)) dz\lesssim\max_{x\in \cA}  \int_{2\D} \left(\frac{\eps}{|z-x|}\wedge 1\right)^qdz.
	\]As the last expression is bounded by $C\eps^{q\wedge 2}$ and $\lim_{\eps\to 0}|\log \eps|\eps^{q\wedge 2}=0$, we  obtain Lemma~\ref{lem:harmonic-circle}.
\end{proof}

The following two crucial lemmas (which are the analogues of Lemmas \ref{Lemma:variance} and \ref{Lemma:covariance} respectively) will be proved in Sections~\ref{sec:finite-variance} and \ref{sec:finite-covariance}.
\begin{lemma}\label{lem:finite-variance}
	There exists a constant $C$ and a sub-polynomial function $\psi$ such that
	\begin{equation}\label{eq:finite-variance}
		\Var[g_\eps(z)]\le \psi(\eps) \quad \forall z\in \D.
	\end{equation}
\end{lemma}

\begin{lemma}\label{lem:finite-covariance}
	Fix a sub-polynomial function $\psi$. We have
	\begin{equation}\label{eq:finite-covariance}
		\Cov[g_\eps(z_1),g_\eps(z_2)]=o_\eps(1)
	\end{equation}uniformly in $z_1,z_2\in \D$  and $|z_1-z_2|>1/\psi(\eps)$.
\end{lemma}

Given these two lemmas, we quickly verify the proof of Theorem \ref{thm:gen}.

\begin{proof}[Proof of Theorem~\ref{thm:gen} given Lemmas~\ref{lem:finite-variance} and \ref{lem:finite-covariance}]
	For $\rho \in C_0^\infty(\D)$, we have
	\begin{align*}
		\Var[(g_\eps,\rho )]=& \int\int  \Cov[g_\eps(z_1) ,g_\eps(z_2)] \rho(z_1)\rho(z_2) dz_1 dz_2.
	\end{align*}
	Let $\psi_1(\eps)$ be the sub-polynomial function in Lemma~\ref{lem:finite-variance} and $\psi_2(\eps)=\psi_1(\eps)\log(1/\eps)$, which is also sub-polynomial. Combining Lemmas~\ref{lem:finite-variance} and~\ref{lem:finite-covariance} and the Cauchy--Schwarz inequality, we have
	$$  \Cov[g_\eps(z_1) ,g_\eps(z_2)]]\le \1_{\{|z_1-z_2|\le 1/\psi_2(\eps)\}}\psi_1(\eps)+o_\eps(1)\1_{\{|z_1-z_2|\ge 1/\psi_2(\eps)\}}.$$
	Since the $o_\eps(1)$ is uniform in $z_1, z_2$ provided that $|z_1 - z_2| \ge 1/\psi_2(\eps)$, we can integrate the above and obtain
	$\lim_{\eps\to0}  \Var[(g_\eps,\rho )]=0$.
	By Lemma~\ref{lem:harmonic-circle}, $(h^{\har},\rho)$ is measurable w.r.t.  $\mu_h$.
\end{proof}

\subsection{Uniform estimate for logarithmic moments for small balls}\label{subsec:gen}
{Suppose we are in the setting of Theorem~\ref{thm:gen}, As in Section~\ref{sec:GMC}, we fix a ${\bf d} \in (\gamma^2/2,d)$ so that $\gamma < \sqrt{2 {\bf d}}$ and $\cE_{{\bf d}} <\infty$ as defined in~\eqref{eq:energy} is finite.} {We may and will assume without loss of generality that $\sigma(\D) \le 1$.}

Both Lemmas~\ref{lem:finite-variance} and~\ref{lem:finite-covariance} will ultimately follow from the proposition below.
\begin{proposition}\label{thm:gen1}
	{For $x\in \cA$ and $\eps\in (0,\frac12)$, let $\fh_x$ be a zero boundary GFF on $B_{2\eps}(x)$.}	
	Let $\mu_x:=\mu_{\sigma,\fh_x}$ be  the GMC of $h_x$ with parameter $\gamma$ over the reference measure  $\sigma$; see~\eqref{eq:GMC}.
	Then there exists a sub-polynomial function $\psi$ depending on $(\gamma,\bf d,\sigma)$ such that
	\begin{equation}
		\label{eq:cond}
		{\E \left[\left|\log\mu_{x}(B_\eps(x)) \right|^4\right] \le \psi(\eps) \quad \textrm{for each }\eps\in(0,\frac12)  \textrm{ and }x\in \cA.}
	\end{equation}
\end{proposition}

{Note that  the law of $\mu_x ( B_\eps(x))$ depends on $x$ as the measure $\sigma$ is not in general translation invariant.}

Our proof of Proposition~\ref{thm:gen1} relies on the following corollary of Theorem~\ref{C:negmom}.
{\begin{lemma}\label{lem:gen}
		{Set $\ell=\frac{\bar \beta-\gamma^2}{\bar  \beta+2\gamma^2}$ with $\bar \beta=\max\{{\bf d},\sqrt{2{\bf d}}\gamma \}$.}
		In the setting of Proposition~\ref{thm:gen1},   for each $p\in (0,1)$,
		there exists a constant {$C=C(\cE_{\bf d}, p, \bf d, \gamma)$}  such that
		\begin{equation*}
			{\E[\mu_{x}(B_\eps(x))^{-p\ell}]   \le C \eps^{-\gamma^2p\ell/2} \sigma(B_\eps(x))^{-p}},  \quad \textrm{for each } \eps\in (0,\frac12) \textrm{ and }x\in \cA.
		\end{equation*}
	\end{lemma}
	\begin{proof}
		Consider the map $y\mapsto (2\eps)^{-1}(y-x)$ from $B_{2\eps}(x)$ to $\D$.
		Let $\sigma_{x,\eps}$ be the pushforward measure of $\sigma|_{B_\eps(x)}$   under this map; {thus $\sigma_{x, \eps} (\D) = \sigma (B_{2\eps} (x)) \le \sigma(\D) \le 1$}.
		Then $\mu_{x}(B_\eps(x))$ can be written as $(2\eps)^{\gamma^2/2} \int_\D e^{\gamma h(z)} \sigma_{x,\eps}(dz)$, where $h$ is a  zero boundary GFF on $\D$.
		
		Let \(X:=\int_\D e^{\gamma h(z)} \sigma_{x,\eps}(dz)\).  To apply Theorem~\ref{C:negmom}, let
		\[
		\cE_{x,\eps}:=\iint_{\D\times \D} |z-w|^{-\bf d}\sigma_{x,\eps}(dw)\sigma_{x,\eps}(dz)=(2\eps)^{\bf d}\iint_{B_\eps(x)\times B_\eps(x)} |z-w|^{-\bf d}\sigma(dw)\sigma(dz).
		\]
		Since $\cE_{x,\eps}\le (2\eps)^{d}\cE_{\bf d}\le \cE_{\bf d}$ for $\eps\in (0,\frac12)$,
		we can apply Theorem~\ref{C:negmom} with {$(\cE',\sigma, s)=(\cE_{\bf d}, \sigma_{x,\eps},p\ell)$ to bound $\E[X^{-p\ell}]$.}
		Since $\sigma_{x,\eps} (\D)=\sigma(B_\eps(x))$ by the definition of $\sigma_{x,\eps}$,  we obtain Lemma~\ref{lem:gen}.
	\end{proof}}
	
	\begin{proof}[Proof of Proposition~\ref{thm:gen1}]
		For each $a>0$,  we have  $|\log x|^4\le c(x+x^{-a})$ for all $x>0$
		for some constant  $c=c(a)$.
		Note that $\E[\mu_{x} (B_\eps(x))]\lesssim \sigma(B_\eps(x))$.
		Fix $a\in (0,\ell)$. Lemma~\ref{lem:gen} shows that
		$$\E[\mu_x (B_\eps(x))^{- a}] \lesssim  \eps^{ - \gamma^2 a /2 } \sigma(B_\eps(x))^{ - a / \ell}.$$
		Using Property (P), we get {$ \E[|\log \mu_x (B_\eps(x))|^4] \lesssim  \eps^{ - \gamma^2 a /2 - q a / \ell}$ where the implicit constant depends on $a$, $\sigma $ and $\gamma$.}
		Since $a$ can be arbitrarily small, this shows that the left hand side is sub-polynomial, as desired.
	\end{proof}
	\subsection{Proof of Lemma~\ref{lem:finite-variance}: the variance estimate}\label{sec:finite-variance}
	By the Cauchy-Schwarz inequality and \eqref{eq:k}, Lemma~\ref{lem:finite-variance} follows from
	\begin{lemma}\label{lem:k-variance}
		There exist a sub-polynomial function $\psi$  such that
		\begin{equation}\label{eq:k-variance}
			\int\Var[ k_\eps(x)] \omega(z,dx)\le\psi(\eps), \quad \forall\, z\in\D.
		\end{equation}
	\end{lemma}
	\begin{proof}
		We will recursive decomposition and idea as in the proof of Lemma~\ref{Lemma:variance} (see Section \ref{sec:proof-Lemma}). Without loss of generality we can assume $ \eps=2^{-n-2}$. Fix $x\in\frac12\D$. We then apply the Markov property of $h$ inductively to decompose $ h $ into
		$$ h= \sum_{i=0}^{n-1} h^{\har}_i+h^{n-1}$$
		where
		$  h^{\har}_i $ ($0\le i\le n-1$) is
		harmonic on $ B_{2^{-i-2}} (x) $, and supported on $ B_{2^{-i-1}} (x) $ (except for $i = 0$); while
		$ h^{n-1}$ is a Dirichlet GFF in $ B_{2^{-n-1}}(x)=B_{2\eps}(x)$. Furthermore, $h^{n-1}$ is independent from $(h^{\har}_i)_{0 \le i \le n -1}$. As in the proof of Lemma \ref{Lemma:variance} (see again Section \ref{sec:proof-Lemma}), applying the intermediate value theorem in each successive ball and the mean value property of $h^{\har}_i$ for each $0 \le i \le n -1$, we can write
		\begin{align*}
			k_\eps(x) &=\sum\limits_{i=0}^{n-1}( h^{\har}_i(\xi_i)-h^{\har}_i(x))+\gamma^{-1} \log \mu_{h^{n-1}}(B_\eps(x))-h_{\eps}^{n-1}(x),
		\end{align*}
		where $ \xi_i\in B_{\eps_{i}} (x)$ comes from the intermediate value theorem, $\mu_{h^{n-1}}=\mu_{\sigma,h^{n-1}}$ and $h^{n-1}_{\eps}$ is the $\eps$-circle average of $h^{n-1}$.
		
		Let
		\begin{equation}
			\label{eq:R_n}
			\Delta_i(x)= h^{\har}_i(\xi_i)-h^{\har}_i(x) \quad\textrm{and}\quad R_{\eps} (x)= \log \mu_{h^{n-1}}(B_\eps(x)).
		\end{equation}
		Then \begin{equation}\label{eq:k-decomposition}
			k_\eps(x)= \sum\limits_{i=1}^{n-1}\Delta_i(x)-h_{\eps}^{n-1}(x)+\gamma^{-1}R_{\eps} (x).
		\end{equation}
		Let $\overline{\Delta}_i(x)= \max_{y_1,y_2\in B_\eps(x)}\{ h^{\har}_i(y_1)-h^{\har}_i(y_2)\}$. By Lemma \ref{Lemma:modulus of harmonic extension},
		\begin{align}\label{eq:variance-individual}
			\E[\Delta^2_i(x)] \le \E [\overline\Delta_i(x)^2]
			\le C 4^{i - n}.
		\end{align}
		Using  \eqref{eq:variance-individual} and the Cauchy--Schwarz inequality, we see that
		\begin{equation}\label{eq:log-variance}
			\var\left[\sum_{i=1}^n \Delta_i(x)\right] \le n \sum_{i=1}^n \E[\Delta_i^2(x)]\le C |\log\eps|.
		\end{equation}
		Observe also that
		$\Var [h^{n-1}_{\eps}(x)] =O(1)$.
		Hence by Cauchy--Schwarz, to finish the proof of Lemma~\ref{lem:finite-variance}, it only remains to find a sub-polynomial function $\psi$ such that
		\begin{equation}\label{eq:key-equation}
			\int \Var \left[ R_{\eps} (x)\right]\omega(z,dx)  \le\psi(\eps).
		\end{equation}
		{From Proposition~\ref{thm:gen1}, we can actually find a sub-polynomial function $\psi$ such that  $\E[R_{\eps} (x)^4]\le \psi(\eps)$ uniform in $x\in \cA$. This concludes the proof of Lemma~\ref{lem:finite-variance}.}
	\end{proof}

	\subsection{Proof of Lemma~\ref{lem:finite-covariance}: the covariance estimate}\label{sec:finite-covariance}
	$\Omega_{\eps,r}=\{(x_1,x_2):|x_1-x_2|\ge\eps^{\frac12},x_1,x_2\in \frac12\D\}$. Recalling \eqref{eq:def-k}, write
	\begin{equation*}
		\E\left[\Cov[g_\eps(z_1),g_\eps(z_2)]\right]=M(z_1,z_2)+R(z_1,z_2),
	\end{equation*}where
	\begin{align}
		& M(z_1,z_2)=\E\left[\iint_{\Omega_{\eps,r}}\Cov\left[ k_\eps(x_1),k_\eps(x_2)\right] \omega(z_1,dx_1)\omega(z_2,dx_2) \right],\nonumber\\
		& R(z_1,z_2)=\E\left[\iint_{|x_1-x_2|<\eps^{\frac12} }\Cov\left[ k_\eps(x_1),k_\eps(x_2)\right] \omega(z_1,dx_1)\omega(z_2,dx_2) \right].\label{eq:def-R}
	\end{align}
	We claim that:
	\begin{lemma}\label{lem:M}
		{We have \(M(z_1,z_2)=o_\eps(1)\) uniformly in $z_1,z_2\in\D$.}
	\end{lemma}
	\begin{lemma}\label{lem:x1-x2-near}
		Fix a sub-polynomial function $\psi$.We have
		\begin{equation}\label{eq:x1-x2-near}
			R(z_1,z_2)=o_\eps(1)
		\end{equation}uniformly in $z_1,z_2\in\D$ with $|z_1-z_2|\ge 1/\psi(\eps)$.
	\end{lemma}
	Note that  Lemma~\ref{lem:M} and Lemma~\ref{lem:x1-x2-near} together will imply  Lemma~\ref{lem:finite-covariance}.
	\begin{proof}[Proof of Lemma~\ref{lem:M}]
		For $x_1,x_2\in \Omega_{\eps,r}$, let $L$ be the straight line orthogonally bisecting $\overline{x_1x_2}$ as in the proof of Lemma~\ref{Lemma:covariance}. Let $ U_i$ be the connected component of $ \D\setminus L$
		containing $ x_i $, $ i=1,2 $. By the domain Markov property, write $h = h^{\har} + h^1 + h^2$, where $h^i$ are independent Gaussian free fields with Dirichlet boundary conditions in $U_i$ ($i = 1, 2$) and $h^{\har}$ is independent, and harmonic in $U_1 \cup U_2$.
		Let $\mu_{h_i}$ ($i=1,2$) be the {GMC measure $\mu_{\sigma,h_i}$} in $U_i$ and
		\begin{align}
			\Gamma_\eps(x_i)&=\gamma^{-1} \log \mu_{h_i}(B_\eps(x_i)) - h_\eps^i(x_i),\nonumber \\
			\Delta_i&=\gamma^{-1} \log \int_{B_\eps(x_i)} e^{\gamma h^{\har}(x)}d\mu_{h_i}(x)- \gamma^{-1} \log \mu_{h_i}(B_\eps(x_i))  - h^{\har}(x_i).\label{eq:Delta}
		\end{align}
		As in the proof of Lemma~\ref{Lemma:covariance} (see Section \ref{sec:cov}),
		\begin{equation}\label{eq:kernel}
			k_\eps(x_i) =  \Gamma_\eps(x_i) +\Delta_i.
		\end{equation} 	
		
		By the independence of $h_1,h_2$, for $|x_1-x_2|>\eps^{\frac12}$,
		\[
		\Cov\left[ \Gamma_\eps(x_1),\Gamma_\eps(x_2)\right] =0.
		\]
		Therefore by \eqref{eq:kernel}
		$$M(z_1,z_2)=\rI(z_1,z_2)+ \rII(z_1,z_2)+\rIII(z_1,z_2)$$
		where
		\begin{align*}
			\rI(z_1,z_2)&=\iint_{\Omega_{\eps,r}}\Cov\left[ \Gamma_\eps(x_1) ,\Delta_2\right] \omega(z_1,dx_1)\omega(z_2,dx_2) ,\\
			\rII(z_1,z_2)&=\iint_{\Omega_{\eps,r}}\Cov\left[ \Delta_1 ,\Gamma_\eps(x_2)\right] \omega(z_1,dx_1)\omega(z_2,dx_2),\\
			\rIII(z_1,z_2)&=\iint_{\Omega_{\eps,r}}\Cov\left[\Delta_1 ,\Delta_2\right] \omega(z_1,dx_1)\omega(z_2,dx_2).
		\end{align*}
		Recall \eqref{eq:Delta}. The intermediate value theorem (applied to $ h^{\har}$) shows that $\Delta_i$ can be rewritten as $h^{\har} (\xi_i) - h^{\har} (x_i)$ for some point $\xi_i \in B_\eps(x_i)$. Hence
		using Lemma \ref{Lemma:modulus of harmonic extension} we obtain for $x_1, x_2 \in \Omega_{\eps, r}$,
		\begin{equation}\label{eq:modulus-bound}
			\Var[\Delta_i]\le \E(\Delta_i^2) \le C\eps^2/|x_2-x_1|^2 \le C \eps,
		\end{equation}
		where $C$ does not depend on $z_1,z_2$. By the Cauchy--Schwarz inequality, we have  $\rIII(z_1,z_2)\le C\eps$.
		
		We now switch out attention to $\rI(z_1,z_2)$ and $\rII(z_1,z_2)$. By the Cauchy--Schwarz inequality (twice),
		\begin{align*}
			\rI(z_1,z_2)&\le \iint_{\Omega_{\eps,r}}\Var^{1/2} [ \Gamma_\eps(x_1)] \Var^{1/2}[ \Delta_2] \omega(z_1,dx_1)\omega(z_2,dx_2)\\
			&\le C\eps^{\frac12} \iint_{\Omega_{\eps,r}}\Var^{1/2} [ \Gamma_\eps(x_1)]  \omega(z_1,dx_1)\omega(z_2,dx_2)\\
			&\le C\eps^{\frac12} \left(\iint_{\Omega_{\eps,r}}\Var [ \Gamma_\eps(x_1)]  \omega(z_1,dx_1)\omega(z_2,dx_2)\right)^{1/2}.
		\end{align*}
		Furthermore, by \eqref{eq:kernel} and \eqref{eq:modulus-bound},
		\begin{align*}
			&\iint_{\Omega_{\eps,r}} \Var[ \Gamma_\eps(x_1)] \omega(z_1,dx_1)\omega(z_2,dx_2)\\
			\le& 2\int\Var[ k_\eps(x_1)] \omega(z_1,dx_1) +2
			\iint_{\Omega_{\eps,r}} \Var[ \Delta_i ] \omega(z_1,dx_1)\omega(z_2,dx_2)\\
			\le & 2\int\Var[ k_\eps(x_1)] \omega(z_1,dx_1)+2C\eps^{\frac12},
		\end{align*}where $C$ is the constant in  \eqref{eq:modulus-bound}.
		A similar estimate holds for $\rII(z_1,z_2)$.
		Now  Lemma~\ref{lem:M} follows from Lemma~\ref{lem:k-variance}.
	\end{proof}
	\begin{proof}[Proof of Lemma~\ref{lem:x1-x2-near}]
		Let
		\begin{align*}
			\rI'(z_1,z_2)&=\iint\1_{\{|x_1-x_2|< \eps^{1/2}\}} \omega(z_1,dx_1)\omega(z_2,dx_2) ,\\
			\rII'(z_1,z_2)&=\iint_{|x_1-x_2|< \eps^{1/2}}\Var\left[ k_\eps(x_1)\right] \omega(z_1,dx_1)\omega(z_2,dx_2) ,\\
			\rIII'(z_1,z_2)&=\iint_{|x_1-x_2|< \eps^{1/2}}\Var\left[ k_\eps(x_2)\right] \omega(z_1,dx_1)\omega(z_2,dx_2).
		\end{align*}
		By the Cauchy-Schwarz inequality,
		\begin{align}\label{eq:R<I}
			R^2(z_1,z_2)
			\le& \left( \iint_{|x_1-x_2|< \eps^{1/2}}\Var^{1/2}\left[ k_\eps(x_1)\right] \Var^{1/2}\left[ k_\eps(x_2)\right]\omega(z_1,dx_1)\omega(z_2,dx_2) \right)^2\nonumber\\
			\le& \rII'(z_1,z_2)\rIII'(z_1,z_2).
		\end{align}
		
		{We first estimate $\rI'(z_1,z_2)$. Assume $\eps$ is small enough such that $\psi(\eps)>3\eps^{1/2}$.
			Then $|z_1-z_2|>\psi(\eps)$ and $|x_1-x_2|<\eps^{1/2}$ imply that either $|z_1-x_2|$ or $|z_2-x_1|$ is bigger than $\psi(\eps)/3$. Note that
			\[
			\iint\1_{\{|x_1-x_2|< \eps^{1/2}, |z_1-x_2|>\psi(\eps)/3\}} \omega(z_1,dx_1)\omega(z_2,dx_2)\le \int \1_{\{|z_1-x_2|>\psi(\eps)/3\}} \omega(z,B_{\eps^{1/2}}(x_2))\omega(z_2,dx_2)
			\]
			Since $\cA$ satisfy property (Q) with exponent $q$ as in~\eqref{eq:propertyQ},
			\[
			\int  \1_{\{|z_1-x_2|>\psi(\eps)/3\}} \omega(z,B_{\eps^{1/2}}(x_2))\omega(z_2,dx_2) \lesssim\eps^{q/2} \psi(\eps)^{-q}= o(\eps^{q/3}).
			\]
			We get a similar estimate if $(z_1,x_2)$ is switched with $(z_2,x_1)$. Therefore	 $\rI'(z_1,z_2)=o(\eps^{\frac q3})$.}

		We now estimate $\rII'(z_1,z_2)$ and $\rIII'(z_1,z_2)$.
		Recall $R_{\eps}(x)$ from the recursive decomposition in \eqref{eq:R_n}. By \eqref{eq:k-decomposition} and \eqref{eq:log-variance}, we have
		\[
		\Var[k_\eps(x)-\gamma^{-1}R_{\eps}(x)]\le C|\log\eps|.
		\]
		Then
		\begin{equation}
			\label{eq:II'<I'}
			\rII'(z_1,z_2)\le C |\log\eps|\rI'(z_1,z_2)+C\iint_{|x_1-x_2|< \eps^{1/2}}\Var\left[ R_{\eps}(x_1)\right] \omega(z_1,dx_1)\omega(z_2,dx_2).
		\end{equation}
		For the second term on the right hand side of \eqref{eq:II'<I'},
		by Proposition~\ref{thm:gen1} and Jensen's inequality $\Var [ R_\eps(x_1)]^2\le \E[R_\eps(x)^4]$.
		Hence this second term is bounded by $\sqrt{\E[R_\eps(x)^4]} I'(z_1, z_2)$.
		Therefore
		$$
		II' (z_1, z_2) \le C\left(|\log\eps| +  \sqrt{\E[R_\eps(x)^4]} \right)  I' (z_1, z_2).
		$$
		As argued at the end of the proof of Lemma~\ref{lem:k-variance},  $\E[R_\eps(x)^4]$ is sub-polynomial.
		Combining with our earlier bound on $I'(z_1, z_2)$,
		we deduce that $\rII'(z_1,z_2)=o_\eps(1)$.
		
		Similarly, $\rIII'(z_1,z_2)=o_\eps(1)$.  Now Lemma~\ref{lem:x1-x2-near} follows from \eqref{eq:R<I}.\qedhere
	\end{proof}

\bibliographystyle{plain}

\end{document}